\documentclass[11pt]{amsart}
\usepackage{amssymb}
\usepackage{eucal}

\begin{document}
\baselineskip=15.5pt

\newcommand{\thmref}[1]{Theorem~\ref{#1}}
\newcommand{\secref}[1]{\S~\ref{#1}}
\newcommand{\defref}[1]{Definition~\ref{#1}}
\newcommand{\lemref}[1]{Lemma~\ref{#1}}
\newcommand{\propref}[1]{Proposition~\ref{#1}}
\newcommand{\corref}[1]{Corollary~\ref{#1}}
\newcommand{\remref}[1]{Remark~\ref{#1}}

\newcommand{\nc}{\newcommand}

\nc{\mc}{\mathcal}
\nc{\on}{\operatorname}
\nc{\Z}{{\mbb Z}}
\nc{\C}{{\mbb C}}
\nc{\Oo}{{\mc O}}
\nc{\D}{{\mc D}}
\nc{\E}{{\mc E}}
\nc{\Ef}{{\mf E}}
\nc{\Ll}{{\mc L}}
\nc{\M}{{\mc M}}
\nc{\Mf}{{\mf M}}
\nc{\Xf}{{\mf X}}
\nc{\N}{{\mc N}}
\nc{\la}{\lambda}
\nc{\V}{\mc{V}}
\nc{\inv}{^{-1}}
\nc{\ol}{\overline}
\nc{\wt}{\widetilde}
\nc{\wh}{\widehat}
\nc{\mb}{\mathbf}
\nc{\mf}{\mathfrak}
\nc{\mbb}{\mathbb}
\nc{\Cx}{{\mbb C}^\times}
\nc{\un}{\underline}
\nc{\sm}{\setminus}
\nc{\boxt}{\boxtimes}
\nc{\ot}{\otimes}
\nc{\Pp}{{\mbb P}}
\nc{\al}{\alpha}
\nc{\pa}{\partial}
\nc{\Ev}{E^{\vee}}
\nc{\No}{{\mc M}^{\on{Id}}}
\nc{\s}{{\mf s}}
\nc{\sbar}{{\ol{\s}}}
\nc{\T}{{\mc T}}
\nc{\Ac}{{\mc A}}
\nc{\XX}{X\times X}
\nc{\K}{\mc K}
\nc{\Sug}{{\mbb S}}
\nc{\Pc}{{\mc P}}
\nc{\Pg}{\on{Pic}^{g-1}}
\nc{\Jac}{\on{Jac}}
\nc{\Ho}{{\rm H}^0}
\nc{\Hi}{{\rm H}^1}
\nc{\Proj}{{\mc Proj}}
\nc{\Conn}{{\mc Conn}}
\nc{\Kern}{{\mc Kern}}
\nc{\undE}{\un{E}}
\nc{\excon}{{\mc E}x\Conn}
\nc{\exhiggs}{{\mc E}x{\mc H}iggs}
\nc{\Oxy}{{\Oo({\mf y}-{\mf x})}}
\nc{\arrtheta}{\stackrel{\rightarrow}{\theta}}
\nc{\gchoosetwo}{\left( \begin{array}{c}g\\2 \end{array}\right)}
\nc{\bra}{\langle}
\nc{\ket}{\rangle}
\nc{\dlog}{\on{d}\log}
\nc{\vphi}{\varphi}
\nc{\Ohalf}{{\Omega_X^{\frac{1}{2}}}}
\nc{\Ominus}{{\Omega_X^{-\frac{1}{2}}}}
\nc{\Thetasing}{\Theta_X^{\text{sing}}}

\title[Theta functions and Szeg\"o kernels]{Theta functions and
Szeg\"o kernels}

\author{David Ben-Zvi}

\address{Department of Mathematics, The University of Chicago,
5734 University Avenue, Chicago, Illinois 60637, USA}

\email{benzvi@math.uchicago.edu}
 
\author{Indranil Biswas}

\address{School of Mathematics, Tata Institute of Fundamental
Research, Homi Bhabha Road, Bombay 400005, India}

\email{indranil@math.tifr.res.in}

\date{}

\begin{abstract}
We study relations between two fundamental constructions associated to
vector bundles on a smooth complex projective curve: the theta
function (a section of a line bundle on the moduli
space of vector bundles)
and the Szeg\"o kernel (a section of a vector bundle on
the square of the curve). Two types of relations are
demonstrated. First, we establish a higher--rank version of the prime
form, describing the pullback of determinant line bundles by
difference maps, and show the theta function pulls back to the
determinant of the Szeg\"o kernel. Next, we prove that the expansion
of the Szeg\"o kernel at the diagonal gives the logarithmic derivative
of the theta function over the moduli space of bundles for a fixed, or
moving, curve. In particular, we recover the identification of the
space of connections on the theta line bundle with moduli space of flat
vector bundles, when the curve is fixed. When the curve varies, we
identify this space of connections with the moduli space of
{\em extended connections}, which we introduce.

\end{abstract}
\maketitle

\section{Introduction}
Let $X$ be a connected smooth projective curve over $\C$ of genus
$g$. There are two fundamental constructions associated to vector
bundles over $X$ of rank $n$ and Euler characteristic zero (hence
degree $n(g-1)$). The first is the nonabelian theta function, which is
a canonical section of the dual of the {\em determinant line bundle}
on the corresponding moduli space. It is nonzero precisely for vector
bundles $E$ which have no sections, that is $\Ho(X,E)=0=\Hi(X,E)$,
which form the complement of the canonical theta divisor $\Theta_X$.
The second starts with such a vector bundle $E$ with no sections, and
constructs a canonical ``kernel function'' on $\XX$, the Szeg\"o
kernel of $E$. This kernel is the unique section of the vector bundle
$E \boxtimes (E^*\ot\Omega_X)$ over $X\times X$ ($\Omega_X$ is the
holomorphic cotangent bundle
of $X$), with only a first order pole along the
diagonal and residue the identity endomorphism of $E$. The Szeg\"o
kernel in the case of line bundles has been extensively studied (see
for example \cite{Hawley, Fay, Raina, KNTY}), as part of the study of
abelian theta functions. In the higher rank case, it has appeared in
the pioneering works of Fay \cite{Fay Szego},
\cite{Fay nonabelian} and
Takhtajan--Zograf \cite{Takh} from an analytic point of view,
where it is related to determinants of Laplace operators and analytic
torsion. (Complex analytically, the role of the Szeg\"o kernel is that
of ``reproducing kernel'' for the holomorphic sections of $E$, or
kernel for the inverse of the Dolbeault operator on $E$.) Our aim in
this paper is to study the Szeg\"o kernel algebraically, investigating
some of its close connections with the theta function, thereby
clarifying some relations between line bundles on moduli spaces and
geometric objects on the curve itself.

\subsection{The nonabelian prime form.}
We will focus on two different mechanisms linking our
protagonists. The first and most direct link is given by difference
maps. To a bundle $E$ is naturally associated the family $E(y-x)$ of
vector bundles on $X$, parametrized by $(x,y)\in\XX$. This family is
classified by a map $\delta_E$ from $\XX$ to the moduli space of
bundles. The pullback of the theta function by $\delta_E$ is a
canonical section of the dual determinant line bundle of this family.
In \thmref{det theorem}, we describe this line bundle canonically in
terms of $E$. The identification, the ``nonabelian prime form'', holds
for arbitrary families of curves and bundles on them (that is, over the
moduli stack of curves and bundles, \corref{stacky}).

In the case of line bundles, this identification is given by the
classical Klein prime form. Moreover, up to a scalar and
multiplication by this prime form, the abelian Szeg\"o kernel and the
pullback of the theta function are in fact equal. (This is sometimes
used as the {\em definition} of the Szeg\"o kernel.) In \thmref{kernel
and theta} we provide a cohomological description of the nonabelian
Szeg\"o kernel and of the nonabelian prime form (for $E$ with no
sections). This makes it easy to compare to the (cohomologically
defined) theta function: the prime form identifies the pullback of the
theta function with the {\em determinant} of the Szeg\"o kernel.  (See
also \remref{prior} for related work.)

\subsection{Twisted cotangent bundles.}
The second mechanism to relate the theta function and Szeg\"o kernel
is through a geometric description of the differential of the theta
function, in terms of connection operators on the curve. More
precisely, we identify the twisted forms of the cotangent bundle which
carry the logarithmic differential of the theta function with spaces
of kernel functions, as well as the corresponding sections.
This is done by describing the behavior of the Szeg\"o kernel along the
generalized Theta divisor.

Let $\Ll$ denote a line bundle on a complex manifold $M$, and
$\Conn_M(\Ll)$ the sheaf of holomorphic connections on $\Ll$. Since
the difference between any two (locally defined) connections on $\Ll$
is a holomorphic one--form, $\Conn_M(\Ll)$ forms an affine bundle
(torsor) for the cotangent bundle $\Omega_M^1$. Such an affine bundle,
equipped with a compatible symplectic form (which for $\Conn_M(\Ll)$
arises as the curvature of a tautological connection), is known as a
{\em twisted cotangent bundle} for $M$ (see \cite{Jantzen}). A
meromorphic section $s$ of $\Ll$ provides a meromorphic section of
$\Conn_M(\Ll)$, which is the coordinate--free form of the logarithmic
differential $\dlog s$. Twisted cotangent bundles on moduli spaces of
curves and of bundles have been studied in \cite{Tyu, Fa, BK, BS}.
It is
known that the twisted cotangent bundle on $\Mf_X(n)$ associated to
the theta line bundle (which we denote $\Conn_X(\Theta)$) is
identified with the moduli space $\Conn_X(n)$ of vector bundles with
flat connection (with its natural symplectic structure).

In \propref{normalized Szego} we describe the sheaf of kernel
functions of which the Szeg\"o kernel is a section, and hence the
behavior of the Szeg\"o kernel along the theta divisor. In
\thmref{dlog theta} we use this description to calculate the
logarithmic derivative of the theta function and provide a new
construction of the isomorphism $\Conn_X(\Theta)\cong
\Conn_X(n)$. Specifically, the restriction of the Szeg\"o kernel to
$2\Delta$ defines a canonical flat connection for every bundle $E$ off
of the theta divisor. The construction develops poles when $E$ hits
$\Theta_X$, which match up precisely with the poles of $\dlog
\theta$. It follows that the restriction of the Szeg\"o kernel is
identified with $\dlog \theta$. (This may be considered an algebraic
analog of the calculation of the logarithmic derivative of the
determinant of the Laplacian in \cite{Fay}, \cite{Fay nonabelian},
\cite{Fay Szego}
and \cite{Takh}.) It also follows that the Szeg\"o kernel defines a
Lagrangian family of flat connections.  (In \secref{identity}, we work
out an identity involving composition of Szeg\"o kernels, which while
not used elsewhere should have independent interest.)

Finally, in \secref{exconns and Szego} we let the curve vary as well
as the bundle, describing twisted cotangent bundles over the moduli
space $\Mf_g(n)$ of curves and bundles on them. To do so we introduce
the notion of {\em extended connections}, which form a twisted
cotangent bundle $\excon_g(n)$ over the moduli space of curves and
bundles. The restriction of the Szeg\"o kernel to $3\Delta$ gives rise
to a canonical extended connection associated to bundles $E$ off of
the theta divisor. In \thmref{dlog theta over curves}, we give a
(unique) identification of the twisted cotangent bundle
$\Conn_g(\Theta)$ over this larger moduli space with $\excon_g(n)$
(due in a different form to \cite{BS}), and identify the logarithmic
differential of the theta function with the restriction of the Szeg\"o
kernel.

\subsection{Organization of the paper.} In Sections \secref{background}
and \secref{flat connections} we provide some background information on
theta functions, twisted cotangent bundles, kernel functions and
connections. In Section \secref{prime section} we discuss the
nonabelian prime form, the difference map, the cohomological
description of the Szeg\"o kernel, and the identification of its
determinant with the theta function. In Section \secref{conns and
Szego} we study the relation of the Szeg\"o kernel to the moduli space
of bundles with connections, by identifying its behavior near the theta
divisor. Finally in Section \secref{exconns and Szego} we introduce
extended connections and extend the results of Section \secref{conns
and Szego} over the moduli space of curves and bundles.

Section \secref{prime section} and Sections \secref{conns and Szego},
\secref{exconns and Szego} can be read independently of each other.

\subsection{Further aspects.}
In \cite{part 2} we study different aspects of extended connections
and their relation with projective structures on $X$. In particular,
we introduce a quadratic map from extended connections to projective
structures, which is a deformation of the quadratic part of the
Hitchin map. Applying this map to the Szeg\"o kernels, we map moduli
spaces of bundles rationally to spaces of projective structures. We
study these maps in detail for line bundles, relating them to several
classical constructions and proving finiteness results. Extended
connections and the quadratic map appear in \cite{Sugawara} from
classical limits of the heat operators for nonabelian theta functions.

\subsection{Remark.}\label{prior} 
After this paper was completed in early 2001, we became aware of the
papers \cite{GG} and \cite{Po}, which prove addition formulae for
non-abelian theta functions, giving a different approach to the
results of \secref{prime section} relating the Szeg\"o kernel and
theta function via the difference map.

\subsection{Notation.}
Throughout this paper, $X$ will denote a smooth connected complex
projective curve (a compact connected Riemann surface) of genus at
least two. All constructions will be algebraic over $\C$. The diagonal
in $X\times X$ consisting of points of the form $(x\, ,x)$ will be
denoted by $\Delta$. For vector bundles $V,W$ on $X$, we will denote
by $V\boxtimes W$ the vector bundle $p_1^* V \bigotimes p_2^* W$ over
$X\times X$, where $p_i$ is the projection to the $i$th factor. For a
coherent sheaf $F$ on $X\times X$, we will use the notation
$F(n\Delta)$ for the sheaf $F\ot \Oo_{X\times X}(n\Delta)$ of
meromorphic sections with only $n$th order pole at the diagonal. The
canonical line bundle (holomorphic cotangent bundle) of $X$ will be
denoted by $\Omega_X$. For a holomorphic vector bundle $E$ over $X$,
its Serre dual, namely $E^*\otimes\Omega_X$, will be denoted by
$E^\vee$. We will often use the same notation for a sheaf and its
space of sections (which is implied should be clear from the context).

\section{Background}\label{background}

\subsection{Theta characteristics.}\label{theta chars}
A {\em theta characteristic} on $X$ is a square root of the canonical
bundle -- that is, a line bundle $\Ohalf$ together with an isomorphism
$(\Ohalf)^{\ot 2}=\Omega_X$. There are $2^{2g}$ possible choices of
$\Ohalf$ on a connected
smooth projective curve of genus $g$. The moduli space
of curves with a choice of theta characteristic is $\Mf_g^{spin}$, the
moduli space of {\em spin curves}.

Tensoring by $\Ominus$ sends bundles of Euler characteristic $0$ to
bundles of degree zero. For a vector bundle $E$, we use the notation
$E_0$ for $E\ot\Ominus$ (so the operation $E\mapsto E_0$ is
well--defined over the moduli space of spin curves.) Note that
$(E_0)^*=(E^{\vee})_0$, so that multiplication by $\Ominus$ converts
Serre duality into sheaf duality.

The ratio of two theta characteristics is a line bundle $\kappa$
equipped with an isomorphism $\kappa^{\ot 2}=\Oo$. Thus $\kappa$
carries a canonical flat connection, inducing the usual connection on
$\Oo$. Using parallel transport for this connection, a canonical
trivialization of $\kappa\boxt \kappa^*$ is obtained on any
neighborhood of $\Delta \subset X\times X$ that contracts to $\Delta$,
in particular, on the infinitesimal neighborhood $n\Delta$ for any
$n$. It follows, for example, that the line bundle
\begin{equation*}
\M_1 \, =\, \Ohalf\boxt\Ohalf(\Delta)
\end{equation*}
on $\XX$ when restricted
to the infinitesimal neighborhood $n\Delta$ for any $n$, is in fact
independent of the choice of theta characteristic. Our use of the
ambiguous notation $\Ohalf$ is thus justified whenever we study
kernel functions of the above form in a neighborhood of
the diagonal. The choice of $\Ohalf$
will appear only when translating statements for Euler
characteristic zero bundles to statements for degree zero bundles.

The tensor power $(\Ohalf)^{\otimes k}$
(for any $k\in\Z$) will be denoted by
$\Omega^{\frac{k}{2}}_X$. Note that the line bundle $\Ohalf$ is Serre
self--dual, $\Ohalf=(\Ohalf)^{\vee}$.

\subsection{Generalized theta functions and determinant bundles.}
\label{det bundles}
We recall some facts about theta functions, determinant bundles and
the canonical theta divisor (see for example \cite{Fa, DN, Qu, BGS}).

Let $\pi:\Xf\to S$ be a smooth family of projective curves over a
connected base $S$, with fibers $X_s$ for $s\in S$. To a coherent
sheaf $E$ on $\Xf$ flat over $S$,
one assigns a line bundle $d(E)$ over
the base $S$, the determinant of the cohomology of $E$, $d(E)=\det
R^{\cdot}\pi_*E$ (EGA III,7). More precisely, this line bundle is
calculated as follows. Locally on $S$ there exists a complex
$$\K^{\cdot}=\{\K^0\stackrel{\pa}{\longrightarrow} \K^1\}$$ of finite
rank free $\Oo_S$--modules, with the property that for any coherent
sheaf $\mc G$ on $S$ there is a natural isomorphism
$$R^i\pi_*(E\ot\pi^*{\mc G})\cong {\rm H}^i(\K^{\cdot}\ot {\mc G}).$$
Such $\K^{\cdot}$ is determined uniquely up to unique
quasi-isomorphism. The determinant line bundle $d(E)$ is then defined
as $d(E)\,=\,\on{det}\K^0\ot(\on{det} \K^1)^*.$ The determinant
functor $d$ descends to a well--defined homomorphism from the
$K$--group of coherent sheaves to the Picard group of $S$ (so that
the determinant of the virtual bundle $E-F$ is $d(E)d(F)^*$). For a
point $s\in S$, the fiber of $d(E)$ over $s$ is identified with the
one dimensional vector space
$$
\wedge^{\text{top}}{\rm H}^0(X_s, E\vert_{X_s})\otimes
\wedge^{\text{top}} {\rm H}^1(X_s, E\vert_{X_s})^*\, .
$$

Now suppose that the Euler characteristic of $E\vert_{X_s}$
vanishes for some (hence every) $s$ so that
$\K^0$ and $\K^1$ are of the same rank. Then the determinant
$\on{det}\pa$ gives a canonical section $\theta$ of the dual bundle
$d(E)^*$, known as the theta function of $E$. The theta function
vanishes precisely for $s\in S$ with ${\rm H}^0(X_s, E\vert_{X_s})
\not= \,0$ (hence ${\rm H}^1(X_s,E\vert_{X_s}) \, \not=\, 0$).

\subsubsection{Moduli spaces.}\label{moduli spaces}
Let $\Mf_X(n)$ denote the moduli space of semistable vector bundles
over $X$ of rank $n$ and Euler characteristic $0$. So the degree of
a vector bundle in $\Mf_X(n)$ is $n(g-1)$. It is known that
$\Mf_X(n)$ is an irreducible normal projective variety. Moreover,
the strictly semistable locus is of codimension greater than one
(unless both $n$ and $g$ are $2$). Its complement
$\Mf_X^s(n)$ defined by all stable vector bundles is smooth.
We will repeatedly take advantage of
the principle that functions, line bundles and twisted cotangent
bundles are not affected by subvarieties of codimension at least
$2$. In other words, the pullback defined by the inclusion of
of the complement of a subvariety of codimension at least two
is an isomorphism. Thus, we will be able to ignore the
strictly semistable locus of $\Mf_X(n)$
and assume all vector bundles involved are stable. 

We will always exclude the special case where both $n$ and $g$
are two.
Although there is no universal vector bundle over $X\times\Mf_X(n)$,
it exists locally (in analytic or \'etale topologies) on the stable
locus $\Mf_X^s(n)\, \subset\, \Mf_X(n)$.
Moreover, all automorphisms of a stable vector bundle are
scalars. It follows that there is a universal vector bundle
over $X\times X\times\Mf_X^s(n)$ of rank $n^2$ such that for
any $E\in \Mf_X^s(n)$, the restriction of the vector bundle
to $X\times X\times\{E\}$ is $E\boxtimes E^*$.

Let $\Mf_X(n)_0$ denote the moduli space of semistable vector bundles
of rank $n$ and degree $0$. If we fix a theta characteristic
$\Ohalf$, then there is an isomorphism 
$$\Mf_X(n)\longrightarrow \Mf_X(n)_0,\hskip.3in E\mapsto
E_0=E\ot\Ominus$$ since tensoring by a line bundle preserves
semistability. Thus over the moduli space of spin
curves, the moduli space of degree zero and Euler characteristic zero
bundles are canonically identified. 

The cotangent bundle to $\Mf_X(n)$ (at stable points $E$) is
identified with the space $\Ho(X,\, \Omega_X\ot\on{End}E)$ of
endomorphism--valued one--forms, or {\it Higgs fields}, on $E$.

\subsubsection{Line bundles on moduli.}
The moduli space $\Mf_X(n)$ carries a universal determinant line
bundle. This is seen as follows. If $E$ is a vector bundle over
$X\times S$ for some scheme $S$ and $L$ a line bundle over $S$, then
it is easy to see from the projection formula that
\begin{equation}\label{det-line-twist}
d(E\otimes\pi^*_S L) \, \cong\, d(E)\otimes L^\chi
\end{equation}
where $\pi_S$ is the projection of $X\times S$ to $S$
and $\chi$ is the Euler characteristic of the restriction of
$E$ to the fibers of $\pi_S$ (which is a
locally constant function on $S$). The determinant line bundle over
$\Mf_X(n)$ is constructed using the local universal families; the
independence of the choices of local families is ensured by
\eqref{det-line-twist}. 

The determinant line bundle on $\Mf_X(n)$ has a concrete geometric
description. The subvariety
$$
\Theta_X \, :=\, \{V\in \Mf_X(n)\, \vert\, \Ho(X,\, V)\,
\not=\, 0\} 
$$
is a divisor, the {\it generalized theta divisor}, that gives the
ample generator of ${\rm Pic}(\Mf_X(n))$ \cite{DN}. Note that
for any $E$ in $\Mf_X(n)$, we have $\Ho(X,\, E)\, =\, 0$ if
and only if $\Hi(X,\, E)\, =\, 0$. This condition also guarantees
that $E$ is semistable. Indeed, if a subbundle $F$ of $E$
violates the
semistability condition, then $\text{h}^0(F)- \text{h}^1(F) >0$, thus
contradicting the condition that ${\rm h}^0(E)=0$ (as
${\rm H}^0(X,F) \, \subset\, {\rm H}^0(X,E)$).
The line bundle over $\Mf_X(n)$ defined by the divisor $\Theta_X$ will
also be denoted by $\Theta_X$, and is canonically identified with the
determinant line bundle.

It is known that the smooth locus of the theta
divisor $\Theta_X$ is precisely the subvariety $\Theta_X^\circ$ of
vector bundles $E$
with ${\rm h}^0(E)={\rm h}^1(E)=1$. The singular locus $\Thetasing\,
\subset\,
\Theta_X$ is of codimension at least $2$ and consists
of all vector bundles $E$ with ${\rm h}^0(E)={\rm h}^1(E)>1$.

\subsection{Twisted cotangent bundles}\label{twisted cotangent bundles}
(Our reference for twisted cotangent bundles is \cite{Jantzen}; see
also \cite{Tyu} for $\Omega$--torsors.)

Let $\Ll$ denote a holomorphic
line bundle over a complex manifold $M$, and $\Conn_M(\Ll)$ the
sheaf of holomorphic connections on $\Ll$. The difference
between any two connections is a (scalar--valued) one--form on $M$.
More precisely,
$\Conn_M(\Ll)$ is an affine bundle for the cotangent sheaf $\Omega_M$ of
$M$ (an $\Omega_M$--torsor). Recall that an affine bundle, or torsor,
is the relative version of the notion of
affine space. Namely given a holomorphic vector bundle $V$ over a
variety $M$, an affine bundle for $V$ over $M$ is a morphism $\pi:A\to
M$, which locally admits a section, equipped with a simply transitive
action of the sheaf of sections of $V$ on the sections of
$A$. (We use the terms affine space for a sheaf, affine
bundle and torsor interchangeably.)

It is convenient to describe this torsor
using the Atiyah exact sequence and the bundle of one--jets.
The Atiyah
bundle of $\Ll$ is the sheaf $\Ac_{\Ll}=\D_{\leq 1}(\Ll,\Ll)$ of
differential operators of order at most one acting on sections of
$\Ll$. It sits in an extension
\begin{equation*}
0\to \Oo \to \Ac_\Ll \to \T_M \to 0
\end{equation*}
of the tangent sheaf by the structure sheaf. The dual sequence is an
extension 
\begin{equation*}
0\to \Omega_M \to \Ac_\Ll^* \to \Oo \to 0
\end{equation*}
The bundle of one--jets $J^1\Ll$ of
sections of $\Ll$, which is an extension
\begin{equation*}
0\to \Omega_M\ot\Ll \to J^1\Ll \to \Ll\to 0
\end{equation*}
of $\Ll$ by the sheaf of $\Ll$--valued differentials, is
obtained from the dual to the Atiyah sequence by tensoring by $\Ll$.

The sheaf $\Conn_M(\Ll)$ of connections on
$\Ll$ is naturally identified with the sheaf of
splittings of any of the above three sequences. In particular,
$\Conn_M(\Ll)$ may be identified with the inverse image of the
section $1\in \Gamma(\Oo)$ in $\Ac_\Ll^*$. From this description the
structure of affine space over the sheaf $\Omega_M$ of differentials
is evident.

Consider the projectivization $\Pp\Ac_\Ll^*=\Pp J^1\Ll$, which is a
projective space bundle containing a projective
subbundle $\Pp\Omega_M$. The
complementary affine space bundle is naturally identified with
$\Conn_M(\Ll)$. A nonzero meromorphic section $s$ of the line bundle
$\Ll$ gives rise to a one--jet $J^1s$, which by projectivization gives
rise to a meromorphic section of $\Conn_M(\Ll)\subset \Pp J^1\Ll$,
which has poles whenever $s$ vanishes or has a pole. This is
called the {\em logarithmic differential} and it is denoted by
$\dlog s$. The resulting connection on $\Ll$
has logarithmic singularities along the divisor of $s$, and can
be characterized as the unique meromorphic connection for which
$s$ is a flat section.

\subsubsection{Symplectic structure.}\label{symplectic}
Let
$$
\pi \, :\, \Conn_M(\Ll)\, \longrightarrow\, M
$$
be the projection map. The pullback bundle $\pi^*\Ll$
over $\Conn_M(\Ll)$ has a
tautological connection $\nabla_\Ll$. To see this
first observe that over
$\Conn_M(\Ll)$ we have a tautological section of
$\pi^*\on{Hom}(\Ll,J^1\Ll)$. There is a natural inclusion
of $\pi^*J^1\Ll$ in $J^1\pi^*\Ll$ defined by the pullback
operation. Combining these two remarks, we obtain a
section of $\on{Hom}(\pi^*\Ll,\pi^*J^1\Ll)$. The
connection on $\pi^*\Ll$ is defined by this section.

Let $\omega_\Ll$ be the curvature of ${\nabla}_L$, which is a
holomorphic $2$-form on $\Conn_M(\Ll)$. It is easy to see that
${\omega}_L$ is a symplectic form. Indeed, if we choose a local
trivialization of $\Ll$, the resulting identification $\Conn_M(\Ll)\to
\Omega^1_M$ takes ${\omega}_L$ to the canonical symplectic form on the
cotangent bundle. It follows that the form $\omega_\Ll$ is changed by
the pullback of
$d\alpha$ under the action of a local section $\alpha$ of $\Omega^1_M$
on $\Conn_M(\Ll)$. In particular, the symplectic form is preserved
precisely by the action of {\em closed} one--forms on $\Conn_M(\Ll)$.
It follows that $\Conn_M(\Ll)$ is a {\em twisted cotangent bundle}:

\subsubsection{Definition.} A twisted cotangent bundle on $M$ is an
$\Omega^1_M$--torsor $A$, equipped with a symplectic form $\omega_A$
which transforms under the translation
action of a local section $\alpha$ of
$\Omega^1_M$ as
$$
\alpha:\omega_A\mapsto \omega_A+df^*\alpha
$$
where $f$ denotes the projection of $A$ to $M$.
 
\subsection{Rigidity.}\label{rigidity}
Twisted cotangent bundles on the moduli spaces of curves and bundles
enjoy a strong rigidity property. The automorphisms of a twisted
cotangent bundle on $M$ are given by closed one--forms on $M$, and the
automorphisms of an $\Omega^1_M$--torsor are given by all one--forms
on $M$. However, as we explain below, the moduli spaces in question
carry no one--forms at all. It follows that any two isomorphisms
between $\Omega^1_M$--torsors on these spaces automatically agree. This
will enable us to compare our descriptions of $\Omega^1_M$--torsors on
moduli spaces with those of \cite{Tyu, Fa, BS}.

\subsubsection{Moduli space of curves.}
It is not hard to verify that the moduli space $\Mf_g$ does not carry
any closed one--forms $\alpha$ for $g>2$.
The first cohomology of $\Mf_g$ vanishes (\cite{HL}), so that $\alpha$
must in fact be exact as the periods of a one-form must be zero.
But $\Mf_g$ does not admit any
nonconstant holomorphic functions if $g\,>\,2$. This follows from the
existence of the Satake compactification of $\Mf_g$, which is a
projective variety in which the complement of $\Mf_g$ has codimension
$2$, for $g>2$. Therefore the form $\alpha$ must vanish identically.
(Note that ${\mathcal M}_2$ is affine, and hence it admits nonzero
exact holomorphic one-forms.)

In fact, a stronger statement holds: $\Mf_g$ does not carry
{\em any} nonzero holomorphic one--forms, provided $g\, >\, 2$.
This follows from the fact, explained to us by E. Looijenga in
\cite{Lo}, that all two--forms on $\Mf_g$ vanish (hence all one--forms
must be closed, and therefore zero.)

\subsubsection{Moduli of vector bundles}
For a fixed line bundle $\xi\in \Pg(X)$, let $\Mf_X(n,\xi) \,\subset\,
\Mf_X(n)$ be the moduli space of semistable vector bundles $E$ with
$\det E \cong \xi$. The nonexistence of holomorphic one--forms on
$\Mf_X(n,\xi)$ follows from unirationality of the moduli space.
It also follows immediately from the usual method of computing
cohomologies using the Hecke transformation.

A similar statement holds for the moduli space
of vector bundles \textit{without fixed determinant}
if we vary $X$ as well. Namely, all the one--forms on $\Mf_X(n)$
are pullbacks of one--forms from
$\Pg(X)$, which are in turn identified with $\Ho(X,\Omega_X)$.
Since there are no holomorphic relative forms over the universal curve
$\Xf\to\Mf_g$, it follows that there is no nonzero one--form
on the universal moduli space $\Mf (n)$ over $\Mf_g$.

\section{Connections and kernels}\label{flat connections}

\subsection{Kernel functions.}
Let $V,W$ be two holomorphic vector bundles over $X$, and let
$\on{Diff}^n(V,W)$ be the sheaf of differential operators of order
$n$ from $V$ to $W$. There is an elegant description of
$\on{Diff}^n(V,W)$ in terms of ``integral kernels'' on $X\times X$,
due to Grothendieck and Sato (see \cite{BS, book}). (This is a
coordinate--free reformulation of
the Cauchy integral formula, describing differentiation in terms of
residues.)
Namely, 
\begin{equation}\label{kernel-op}
\on{Diff}^n(V,W)\,\cong\, 
\frac{W\boxtimes V^\vee((n+1)\Delta)}{W\boxtimes
V^\vee}\, \cong\, W\boxtimes (V^*\otimes \Omega_X)
((n+1)\Delta))\big\vert_{(n+1)\Delta}.
\end{equation}
Here we consider $n$th order differential operators as a bimodule over
$\Oo_X$ (via left and right multiplication), and hence as a
coherent torsion sheaf on
$\XX$, supported on $(n+1)\Delta$. In other words, we consider kernel
functions of the form $\psi(z,w) dw$ on $X\times X$, with poles of
order $\leq n+1$ allowed on the diagonal, and quotient out by regular
sections. These act as differential operators by sending $\psi(z,w)dw$
to the operator $\delta_\psi$ defined by
$$
\delta_\psi(f(z)) \,=\, \on{Res}_{z=w} \langle \psi(z,w),
f(z)\rangle dw
$$
where $f(z)$ is a local section of $V$ and $\langle-,-\rangle$
is the contraction of $V^*$ with $V$. So, $\langle \psi(z,w),
f(z)\rangle dw$ is a local section of $W\boxtimes \Omega_X
((n+1)\Delta))$ and hence its residue is a local section of $W$.

\subsubsection{The de Rham kernel}\label{de-rham}
Let $d\,:\, \Oo_X\, \longrightarrow\, \Omega_X$ be the de Rham
differential. Using \eqref{kernel-op}, this operator defines a
holomorphic section $\mu(d)$ of $\Omega_X\boxtimes \Omega_X (2\Delta)$
over $2\Delta$. Since the symbol of the differential operator $d$ is
the section of ${\mc O}_X$ defined by the constant function $1$, the
restriction of $\mu(d)$ to $\Delta$ is the constant function $1$
(considered as a section of $\Omega_X\boxtimes
\Omega_X(2\Delta)|_\Delta={\mc O}_\Delta$, by Poincar\'e
adjunction formula).

More generally, consider the line bundle
$$\M_{\nu} = \Omega_X^{\frac{\nu}{2}}\boxt
\Omega_X^{\frac{\nu}{2}}(\nu\Delta)$$ over $X\times X$, where $\nu\in
\Z$ is an integer. By \secref{theta chars}, it follows that the
restriction $\M_{\nu}|_{n\Delta}$ is independent of the choice of
theta characteristic for any integers $\nu,n$.

We will throughout use the convention that the normal bundle of
$\Delta$ in $X\times X$ is identified with
the tangent bundle $T\Delta$ using the projection
to the first factor. In other words, if $p_1$
denotes the projection to the first factor of $X\times X$, then
the isomorphism is the composition of the isomorphism of the
normal bundle with $p^*_1TX$ defined by the differential of $p_1$
with the natural identification of $T\Delta$
with $TX$. Note that if we use the projection to the
second factor, then the isomorphisms differ by a sign.

Let $\sigma\,:\, X\times X\to X\times X$ be the involution
defined by $(x,y) \mapsto (y,x)$.
Note that the adjunction formula identifies the
restriction of the line bundle ${\mathcal O}_{X\times X}(\Delta)$
to $\Delta$ with the normal bundle $N$ of $\Delta$ in $X\times X$.
Using the ordering of the product $X\times X$, the normal bundle
$N$ gets identified with the tangent bundle of the divisor
$\Delta$. The tangent bundle $T\Delta$ is identified with
$TX$ using the projection $p_2$. On the other hand, projecting
the line subbundle $(p^*_2TX)\vert_\Delta\, \subset\,
(T(X\times X))\vert_\Delta$ to the quotient
$N$ of $(T(X\times X))\vert_\Delta$ we get an isomorphism of
$(p^*_2TX)\vert_\Delta$ with $N$. Note that the isomorphism of $N$
with $T\Delta$ changes by multiplication with $-1$ if the two
factors in $X\times X$ are interchanged (that is, if we use $p_1$
instead of $p_2$).

Consequently, the restriction of $\M_\nu$ to $\Delta$ is identified
with $\Omega^{\otimes \mu}_X\otimes N^{\otimes \mu}$. Using
the above identification of $T\Delta$ with $N$ the line bundle
$\Omega^{\otimes \mu}_X\otimes N^{\otimes \mu}$ is trivialized.
It is easy to see that for any $\nu$ there is a unique
trivialization $\mu_\nu$ of $\M_\nu|_{2\Delta}$ such that
\begin{enumerate}

\item{} $\M_\nu|_{\Delta}\,\cong \,\Oo_X$ (that is,
$\mu_\nu|_\Delta \, =\, 1$);

\item{} the trivialization is symmetric, respecting the
identification $\M_\nu\cong \sigma^*\M_\nu$ (in other words,
$\sigma^*\mu_\nu= (-1)^\nu \mu_\nu$).
\end{enumerate}
The appearance of the factor $(-1)^\nu$ is due to the above
remark that the isomorphism of $N$ with $T\Delta$
changes sign as the two factors are interchanged.
In particular, $\mu_\nu=\mu_1^{\ot \nu}$ and $\mu_2 = \mu(d)$,
where $\mu(d)$ is the section of ${\mc M}_2$ over $2\Delta$
defined by the de Rham differential $d$.

\subsection{Connections.}
Let $E$ denote a vector bundle over $X$ of rank $n$, and $\Conn_X(E)$
the sheaf of holomorphic connections on $E$. Note that $\Conn_X(E)$
will have no global sections unless $E$ has degree zero.

A connection
$\nabla$ on $E$ may be described as a first--order differential
operator $\nabla:E\to E\otimes \Omega_X$ whose symbol is the identity
automorphism of $X$. (This condition on symbols is equivalent to the
Leibniz identity.)
It follows from the differential operators--kernels dictionary in
\eqref{kernel-op} that giving a holomorphic
connection on $E$ is equivalent to giving a section of
$$
\M_2(E):=(E\otimes \Omega_X) \boxt (E^*\ot \Omega_X) (2\Delta)
$$
over $2\Delta$ whose
restriction to $\Delta$ is the identity section of
$$
\M_2(E)|_\Delta\,\cong\, \on{End} E.$$ We also see that the difference
between any two connections is a section of
$\Omega_X\ot\on{End}E$. Thus $\Conn_X(E)$ is an affine space for the
space $\Ho(X,\, \Omega_X\ot\on{End}E)$ of Higgs fields on $E$.

More generally, let
$$\M_\nu(E)=(E\boxt E^*) \otimes{\M}_\nu = (E\otimes
\Omega_X^{\frac{\nu}{2}} ) \boxt (E^*\otimes
\Omega_X^{\frac{\nu}{2}})(\nu\Delta).$$ 
On $2\Delta$, we obtain an extension
$$0\to \on{End}E\ot\Omega_X\to \M_\nu(E)|_{2\Delta}\to \on{End}E\to
0.$$ Since $\M_\nu$ is trivialized on $2\Delta$ (by the section
$\mu_\nu$ in \secref{de-rham}), we may identify the sheaves
$\M_\nu(E)|_{2\Delta}$ above for all integers $\nu $. Consequently,
$\Conn_X(E)$ coincides with the space of sections of
$\M_\nu(E)|_{2\Delta}$ lifting $\on{Id}\in\on{End}E$, and the
difference between any two such sections is naturally a Higgs field.

\subsubsection{Remark.} 
Recall that, in Grothendieck's formulation, a connection $\nabla$ on
$E$ is the data of identifications of the fibers of $E$ at
infinitesimally nearby points. More precisely, a connection on $E$ is
an isomorphism between the two pullbacks $p_1^*E$ and $p_2^*E$ on the
first--order infinitesimal neighborhood $2\Delta$ of the diagonal,
which restricts to the identity endomorphism of $E$ on the
diagonal. Note that since
$$\un{Hom}_{\Oo_{X\times X}}(p_2^*E,p_1^*E)= E\boxt E^*,$$ this is
equivalent to the data of a section, which we denote $k_\nabla$, of
$E\boxt E^*$ on $2\Delta$, whose restriction to $\Delta$ is the
identity. This is the case $\nu=0$ of the above construction.

\subsubsection{The Atiyah bundle.}\label{Atiya-bundl.-def.}
Another related point of view on connections is due to Atiyah,
\cite{At}. For a vector bundle $E$ over $X$, the symbol map sends
$\D_{\leq 1}(E,E)$, the sheaf of differential operators of order $1$,
to ${\rm Hom}(E\ot\Omega_X,E)$.
The \textit{Atiyah bundle} ${\rm At}(E)$ of $E$
consists of differential operators with {\em scalar} image. So
${\rm At}(E)$ is the
inverse image of the scalars $\T\subset {\rm Hom}(E\ot\Omega_X,E)$ in
$\D_{\leq 1}(E,E)$. Thus, $\text{At}(E)$ fits in an exact sequence
$$
0\, \longrightarrow\, \on{End}V\, \longrightarrow\,\text{At}(E)
\, \longrightarrow\, \T_X \, \longrightarrow\, 0\, .
$$
In the language of kernels, the Atiyah bundle appears
as follows:
$$
\rm{At}(E)=\{\psi\in\Gamma(E\boxt E^{\vee}(2\Delta)|_{2\Delta})
\;\big\vert\; \psi|_{\Delta}\in \Omega_X\cdot\on{Id}_E\}\, .
$$
(Here we
use the identification $E\boxt E^{\vee}(2\Delta)|_{\Delta}\cong
\Omega_X\ot \on{End}E$ of symbols.) The sheaf $\rm{At}(E)\ot\Omega_X$
is naturally identified with the subsheaf of
$\M_2(E)|_{2\Delta}$ consisting of all sections whose restriction
to the diagonal is contained in the subsheaf
$\Oo\cdot\rm{Id}\subset \on{End}E$.
Hence again we see that a holomorphic connection on $E$
is simply a splitting of the Atiyah sequence.

\subsection{Twisted connections.}
Let $E$ denote a rank $n$ vector bundle on $X$, and consider the sheaf
$$\M(E)=E\boxt\Ev(\Delta)$$ on $\XX$. For an arbitrary choice of theta
characteristic, there is a canonical identification 
$$\M(E)=\M_1(E_0).$$ 

The affine space $\Conn_X(E_0)$ of connections on $E_0=E\ot\Ominus$
depends only on $E$, not on the choice of $\Ohalf$: it is identified
with (sections of) the affine bundle
$$\Conn_X(E_0)=\{s\in \M(E)|_{2\Delta} \big\vert\,\,
s|_\Delta=\on{Id}_E\}.$$ This can also be seen as follows: if $\kappa$
is a line bundle with a given flat connection, then there is a
canonical isomorphism between $\Conn_X(F)$ and
$\Conn_X(F\ot\kappa)$. Since the ratio of two theta characteristics is
such a $\kappa$ (\secref{theta chars}), we see that the affine bundles
$\Conn_X(E_0)$ are independent of the choice.

We will refer to elements of $\Conn_X(E_0)$
as {\em twisted connections} on $E$. In particular $\Conn_X(E_0)$ will
have no global sections unless $E$ has Euler characteristic zero (so
that $E_0$ has degree zero.)

\subsection{The Szeg\"o kernel.}
There is a canonical kernel function associated to vector bundles $E$
off of the theta divisor: the nonabelian Szeg\"o kernel. (The proof
of the uniqueness and existence below is straightforward --- a
stronger version is given in \propref{kernels}.)

\subsubsection{Definition.}\label{def of Szego} 
Let $E\in \Mf_X(n)\sm\Theta_X$ be a vector bundle with
${\rm h}^0(E)={\rm h}^1(E)=0$. The Szeg\"o kernel of $E$ is
the unique section $\s_E\in \Ho(X\times X,\M(E))$
with $\s_E|_{\Delta}=\on{Id}_E$.

Equivalently, we will also consider the Szeg\"o kernel as
a meromorphic section of $E\boxt\Ev$ with a pole of order
exactly one on the diagonal. The that case, the above
condition $\s_E|_{\Delta}=\on{Id}$ translates into the
condition that the residue of the meromorphic section
is the identity automorphism of $E$.

\section{Nonabelian Prime Forms and
the Szeg\"o kernel}\label{prime section}
\subsection{The Prime Form.}\label{the prime form}
Let $\Jac=\Mf_X(1)_0$ denote the Jacobian variety of $X$. Consider
the difference map
$$
\delta\, :\, \XX\to \Jac ,\hskip.3in (x,y)\mapsto y-x\, .
$$
Its image is a surface in the Jacobian referred to as $X-X$. This map
classifies a natural line bundle on
$X\times(\XX)$. The line bundle on $X\times(\XX)$ in question,
which we will denote by $\Oxy$, is defined as follows:
$$
{\Oxy}\, :=\, \Oo(\Delta_{02}-\Delta_{01})\, .
$$
It may also be characterized by the following conditions:
\begin{enumerate}
\item{} for any point $(x,y)\, \in\, \XX$, the restriction
of $\Oxy$ to $X\times (x,y)$ is isomorphic to the line bundle
${\mathcal O}_X(y-x)$ over $X$;

\item{} for a point $z\in X$,
the restriction of $\Oxy$ to $z\times \XX$ is isomorphic
to ${\mathcal O}_X(-z)\boxtimes{\mathcal O}_X(z)$.
\end{enumerate}
Here and in what follows, we use the following notations for the cube
$X\times\XX$: we label the three copies of $X$ by $0,1,2$ and let
$\Delta_{ij}\subset X\times\XX$ denote the diagonal where the points
labelled $i$ and $j$ collide. Let $p_i\, :\, \XX\, \longrightarrow\,
X$, $i\, =\, 1,2$ denote the projection to the $i$-th factor.

Let $\on{Pic}^d(X)$ denote the Picard variety of
degree $d$ line bundles on $X$.
For any line bundle $\Ll\in\on{Pic}^d(X)$, we
translate the difference map by $\Ll$, giving rise to
$$
\delta_{\Ll}:\XX \to \on{Pic}^d(X)\, .
$$
This map classifies the line bundle $\Ll\ot\Oxy$
on $X\times(\XX)$. The inverse image of $\{\Ll\}$ is the
diagonal $\Delta\subset\XX$.

Consider now the variety $\Pg(X)=\Mf_X(1)$ of line bundles
over $X$ of degree
$g-1$. For $\Ll\in\Pg(X)$, it is well--known (see e.g. \cite{Raina})
that the pullback $\delta_\Ll^*\Theta_X$ of the theta line bundle on
$\Pg(X)$ is isomorphic to $$\M(\Ll)=\Ll\boxt \Ll^{\vee}(\Delta).$$
However, while the restriction $\M(\Ll)|_\Delta$ is canonically
trivial, the restriction of the pullback
$\delta_\Ll^*\Theta_X|_{\Delta}$ is naturally identified with the
trivial line bundle $\Oo\otimes \Theta_X|_{\Ll}$ (since the values of
theta functions at $\Ll$ are not numbers). Hence we must tensor
$\M(\Ll)$ by the complex line $\Theta_X|_\Ll$ to make the
isomorphism canonical:

\subsubsection{Definition.}\label{prime form} 
The unique isomorphism
$$E_\Ll: \M(\Ll)\ot\Theta_X|_{\Ll} \longrightarrow
\delta_{\Ll}^*\Theta_X$$ of line bundles on $\XX$, which restricts to
the identity on the diagonal, is known as the {\em prime form} of $\Ll$.

\subsubsection{Translation to degree zero and the classical prime
form.}
\label{Klein prime} 
Let us pick a canonical basis $A_1,\dots,A_g,B_1,\dots,B_g$ in ${\rm
H}_1(X,\Z)$. This choice determines a theta characteristic $\Ohalf$ on
$X$, known as Riemann's constant. It is characterized by the
property that the $\Ohalf$ translate of the divisor
$\Theta_0$ of Riemann's theta
function on $\Jac$ (defined using the polarization $A_i,B_i$) is the
canonical theta divisor $\Theta_X\subset \Pg (X)$ (the Abel--Jacobi
image of $\on{Sym}^{g-1}X$).

The classical Klein prime form of $X$ (for detailed discussions see
\cite{Fay}, \cite{Tata II}) is a section of 
$$\Ominus\boxt \Ominus,$$ with values in $\delta^*\Oo(\Theta_0),$ and
with divisor the diagonal. Dividing by this section defines an
isomorphism $E:\Ohalf\boxt\Ohalf(\Delta)\to \delta^*\Oo(\Theta_0)$. To
recover the isomorphisms $E_{\Ll}$ of \defref{prime form}, note that
$\delta^* \Theta_0=\delta_\Ominus^* \Theta_X$. Since the ratio between
the pullbacks of $\Theta_X$ by $\delta_\Ll$ and $\delta_\Ohalf$ is
$\Ll_0\boxt (\Ll_0)^*= \Ll_0\boxt \Ll^{\vee}_0,$ dividing by the Klein
prime form may be used in defining the isomorphisms $E_\Ll$ for all
$\Ll$.

\subsection{Abelian Szeg\"o kernels and theta functions.}
Suppose $\Ll\in\Pg(X)\sm\Theta_X$ is a line bundle with
${\rm h}^0(\Ll)={\rm h}^1(\Ll)=0$. By \defref{def of Szego} we have a
canonical section
$$\s_\Ll\in\Ho(\XX, \M(\Ll)),$$
the Szeg\"o kernel of $\Ll$. On the other hand, by \defref{prime form}
the prime form of $\Ll$ provides an isomorphism
$$E_\Ll:\M(\Ll)\ot \Theta_X|_{\Ll}\longrightarrow
\delta_\Ll^*\Theta_X,$$ so the theta function pulls back to a
section of $\M(\Ll)$. Normalizing the pullback to have value
$1$ on the diagonal, we recover the well--known formula for the
Szeg\"o kernel,
\begin{equation}\label{formula for Szego}
\s_\Ll(x,y)=\frac{\theta(y-x+\Ll_0)}{\theta(\Ll_0)E(x,y)}.
\end{equation}
Here we use Riemann's constant $\Ohalf$ (\secref{Klein prime}) to
translate $\theta$ and $\Ll$ to $\Jac$, whose group structure is
written additively. In other words, up to the isomorphism given by
the prime form and multiplication by a scalar (the value
$\theta(\Ll_0)$), the Szeg\"o kernel is the pullback of the theta
function by the difference map of $\Ll$.

The Szeg\"o kernel for the line bundle
${\Omega}^{\frac{1}{2}}_{{\mathbb C}{\mathbb P}^1}$ over
${\mathbb C}{\mathbb P}^1 = {\mathbb C}\cup\{\infty\}$ is
$$
\frac{\sqrt{dz_1}\otimes \sqrt{dz_2}}{z_2-z_1}
$$
where $(z_1,z_2)$ is the coordinate on ${\mathbb C}^2$.

\subsection{The Nonabelian Prime Form.}\label{nonabelian prime}

Fix a holomorphic vector bundle $E$ of rank $n$ over $X$. 
The projection
of $X\times(\XX)$ to $X$ (respectively, $\XX$) will be denoted by $p_X$
(respectively, $\Pi$). 
Let $\undE=p^*_X E$ and
\begin{equation}\label{E Oxy}
\E\, :=\, \undE \otimes {\Oxy} 
\end{equation}
be the vector bundle over $X\times(\XX)$. In other words, $\E$
is the family of vector bundles $E(y-x)$ over $X$ parametrized by
$\XX$. 
Let $d(\E)$ be the determinant line bundle over $\XX$ for
the family $\E$ of vector bundles over $X$. 
(Its dual is the pullback of the theta line bundle on the moduli space
of vector bundles, by the map $\delta_E$ classifying the family $\E$.)

We recall that $d(\E) \,=\, \det R^0\Pi_*\E\bigotimes (\det
R^1\Pi_*\E)^*$, and hence the fiber of $d(\E)$ over any $(x,y)\in \XX$
is ${\bigwedge}^{\rm top}\Ho(X,\, E\ot\Oo(y-x))
\bigotimes\left({\bigwedge}^{\rm top}\Hi(X,\, E\ot\Oo(y-x))\right)^*$.
Let $\det E$ denote the line bundle ${\bigwedge}^n E$ on $X$.

\subsubsection{\bf{Proposition}.}\label{det prop}
The dual determinant bundle $d(\E)^*$ over $\XX$ is 
isomorphic to the line bundle $\det E\boxt \det E^{\vee}(n\Delta)$.

\subsubsection{Introducing parameters.}
We will prove a stronger version of \propref{det prop}, valid for
arbitrary families of curves and vector bundles on them (i.e., over
the moduli stack). Hence in the rest of \secref{nonabelian prime} we
work over a fixed (complex) base scheme (or base analytic space) $S$,
and all constructions are relative to $S$. Thus $\pi:\Xf\to S$ will
denote a smooth projective curve over $S$, of genus $g$ (i.e., a
smooth family of compact Riemann surfaces of genus $g$ parametrized by
$S$). For any $s\, \in\, S$, the curve ${\pi}^{-1}(s)$ is denoted by
$X_s$. In this section all products are fiber products over $S$, so
that the projection
\begin{equation}\label{pi-2}
\pi_2 \,:\, \Xf\times\Xf\,\longrightarrow\, S
\end{equation}
will denote the family
$\pi\times_S\pi:\Xf\times_S\Xf$ of algebraic
surfaces $\{X_s\times X_s\}$ and
$\Delta\subset \Xf\times \Xf$ will denote
the relative diagonal. Let $\Pi\, :\,
\Xf\times \Xf\times \Xf\,\longrightarrow\, \Xf\times \Xf$ be the
projection to the second and third factors. Note that the family
of Riemann surfaces over the fiber product $\Xf\times \Xf$
defined by the projection $\Pi$ is the pullback of the
family $\Xf$ over $S$ by the obvious projection of
$\Xf\times \Xf$ to $S$.

The projection of
$\Xf\times \Xf\times \Xf$ to the first factor will be denoted by
$p_X$. Let $\Delta_{01}$ (respectively,
$\Delta_{02}$) denote the diagonal
divisor in $\Xf\times(\Xf\times \Xf)$ defined by all points whose
first coordinate coincides with the second (respectively, third)
coordinate.

Let $E$ denote a vector bundle of rank $n$ over $\Xf$ (i.e., a
holomorphic family of vector bundles $E_s$ over $X_s$), and $d(E)$ its
determinant line bundle over $S$. By $E^{\vee}$ we denote the
relative Serre dual bundle $E^*\ot\omega_{\pi}$, where $\omega_{\pi}$
is the relative dualizing sheaf (relative canonical line bundle) for
the morphism $\pi:\Xf\to S$.

Set $\undE \,:=\, p_X^*E$
and consider the following analogue of \eqref{E Oxy}
$$
\E \, =\, \undE\otimes \Oo(\Delta_{02}-\Delta_{01})
$$
over $\Xf\times (\Xf\times \Xf)$. The determinant bundle $d(\E)$ for
the morphism $\Pi$ is a line bundle over $\Xf\times \Xf$.

\medskip
\subsubsection{\bf{Theorem}.}\label{det theorem}
The theta (or dual determinant) line bundle $(\Theta_X)_{\E}
\,=\, d(\E)^*$ is
canonically isomorphic to
$$
\det E \boxt\det E^{\vee}(n\Delta) \ot \pi_2^*d(E)^*
$$
(where $d(E)^*=(\Theta_X)_E$ is the determinant line bundle
on $S$, and $\pi_2$, as in \eqref{pi-2}, is the
natural projection of $\Xf\times \Xf$ to $S$).

\subsubsection{Remark: The universal situation.}
The theorem is equivalent to the following universal statement on the
moduli {\em stack} of curves and bundles. (A suitably modified
statement also holds over the moduli {\em space} of stable bundles.) Let
$\un{\Mf}_g(n)$ denote the moduli stack of pairs $(X,E)$ of a smooth
projective curve of genus $g$ with a vector bundle
over it of rank $n$, and
$(\Xf,\Ef)$ the universal curve and bundle. Let
$$\delta_\Ef:\Xf\times_{\un{\Mf}_g}\Xf\times_{\un{\Mf}_g}\un{\Mf}_g(n)
\to\un{\Mf}_g(n)$$ denote the universal difference map, classifying
the vector bundle $\Ef\boxt\Ef^{\vee}$.

\subsubsection{Corollary.}\label{stacky} There is an isomorphism
$$\delta_{\Ef}^*\Theta_X \cong \det\Ef\boxt\det\Ef^\vee(n\Delta)\ot
\pi_2^*\Theta_X$$ of line bundles on the square of the
universal curve over the moduli stack of curves and bundles
$\un{\Mf}_g(n)$.

\subsubsection{The Deligne pairing.}
The proof of \thmref{det theorem} relies on the Deligne pairing. We
summarize here the properties of this construction that we will need,
taken from \cite[pp. 98, 147]{De} and \cite[pp. 367--368]{BM}. Given a
family of curves $\pi: \Xf\to S$ and two line bundles $L$ and $M$ on
$\Xf$, the Deligne pairing assigns a line bundle on the base $S$,
denoted by $\langle L\, , M\rangle$. This pairing is a refined version
of the pushforward of the product of Chern classes: the Chern class
$c_1(\langle L\, , M\rangle)$ is the image, under the Gysin
homomorphism ${\rm H}^4(\Xf ,{\mbb Q})\to {\rm H}^2(S,{\mbb Q})$, of
$c_1(L)c_1(M)$. There are canonical isomorphisms as follows:
\begin{enumerate}
\item[a.] $\bra L,M\ket=\bra M,L\ket$, $\langle L,M\ot
N\rangle=\langle L,M\rangle\ot \langle L,N\rangle$, $\langle L,
M^*\rangle=\langle L, M\rangle^*$, and $\langle L,\Oo_X\rangle =
\Oo_S$.
\item[b.] $\langle L,M\rangle=d(L\ot M)d(\Oo_X)d(L)\inv d(M)\inv$.
\item[c.] For pairs of vector bundles $E_0,E_1$ and $F_0,F_1$ of the
same rank, $$\langle \det(E_0-E_1),\det(F_0-F_1)\rangle =
d((E_0-E_1)\ot(F_0-F_1))$$ (where the differences denote virtual
bundles).
\item[d.] For a relatively positive divisor $D\subset X$ with
structure sheaf $\Oo_D$, there is an isomorphism $\langle
L,\Oo_X(D)\rangle = d(L\ot\Oo_D)\ot d(\Oo_D)\inv$.
\end{enumerate}

\subsubsection{Proof of \thmref{det theorem}.}
There are short exact sequences
$$0\to \Oo(-\Delta_{01})\to \Oo\to \Oo_{\Delta_{01}}\to 0 \hskip.2in
\mbox{and}$$
$$0\to \Oo \to \Oo(\Delta_{02}) \to \omega_{\Delta_{02}}^* \to 0,$$ by
the Poincar\'e adjunction formula. 
Since the determinant line is
multiplicative in short exact sequences, we obtain
\begin{equation}\label{iso-deligne-exact}
d(\Oo(-\Delta_{01}))=d(\Oo),\hskip.3in
d(\Oo(\Delta_{02}))=d(\Oo)\otimes p_2^*\omega^*_\pi
\end{equation}
where $p_i$ is the projection of $\Xf\times \Xf$ to the
$i$th factor. (Note that for a sheaf whose support has relative
dimension zero, the determinant line is simply the
top exterior power of the pushforward.)

We next calculate the Deligne pairing $\bra
\Oo(\Delta_{02}),\Oo(-\Delta_{01}) \ket$ in two different ways. Using
$(d)$ (and symmetry) we obtain
\begin{equation*}
\bra \Oo(\Delta_{02}),\Oo(-\Delta_{01}) \ket =
\Oo(-\Delta),
\end{equation*}
as $d(\Oo_{\Delta_{02}})$ is the trivial line bundle,
while using $(b)$ we obtain
\begin{equation*}
\bra \Oo(\Delta_{02}),\Oo(-\Delta_{01}) \ket =
d(\Oo(\Delta_{02}-\Delta_{01}))\otimes d(\Oo(-\Delta_{02}))^*
\otimes d(\Oo(\Delta_{01}))^*\otimes d(\Oo).
\end{equation*}
Comparing these two expressions for
$\bra \Oo(\Delta_{02}),\Oo(-\Delta_{01}) \ket$
and using \eqref{iso-deligne-exact} we find
\begin{equation}\label{d of oxy}
d(\Oo(\Delta_{02}-\Delta_{01}))^*= p_2^*\omega_\pi\ot d(\Oo)^*\ot
\Oo(\Delta).
\end{equation}

We now calculate $\langle \det\undE,\Oo(\Delta_{02}-\Delta_{01})
\rangle$ in two ways. First, using $(a)$ and $(d)$, we may
write
\begin{eqnarray}\label{deligne-isomorphism}
\nonumber 
\langle \det\undE,\Oo(\Delta_{02}-\Delta_{01})
\rangle &=& 
\bra \det\undE, \Oo(\Delta_{02})\ket \otimes\bra \det\undE,
\Oo(\Delta_{01})\ket^*\\
\nonumber 
&=& d(\det\undE\otimes\Oo_{\Delta_{02}})
\otimes d(\Oo_{\Delta_{02}})^*\otimes
\left(d(\det\undE\otimes\Oo_{\Delta_{01}})
\otimes d(\Oo_{\Delta_{01}})^*\right)^*\\
&=& p_2^*\det E \ot p_1^*\det E^*
\end{eqnarray}
(since $d(\det\undE\otimes\Oo_{\Delta_{0i}})\,=\,
p_i^*\det E$).
On the other hand, using $(a)$ and $(c)$ (and the fact that
$\text{det}$ and determinant are homomorphisms from the Grothendieck
$K$-group $K(\Xf)$ to the Picard group of line bundles over
$\Xf$ and $S$ respectively), we write
\begin{eqnarray*}
\langle \det\undE,\Oo(\Delta_{02}-\Delta_{01})
\rangle &=& 
\bra \det(\undE-\Oo^{\oplus n}),
\det(\Oo(\Delta_{02}-\Delta_{01})-\Oo)\ket \\
&=& d((\undE-\Oo^{\oplus n})\otimes
(\Oo(\Delta_{02}-\Delta_{01})-\Oo))\\
&=& d(\E)\otimes d(\undE)^*\otimes d(\Oo(\Delta_{02}-\Delta_{01}))^{-n}
\otimes d(\Oo)^n .
\end{eqnarray*}
Comparing these and using \eqref{d of oxy} and
\eqref{deligne-isomorphism} we obtain
\begin{eqnarray*}
d(\E)^*&=& p_1^*\det E\ot p_2^*\det E^* \ot (p_2^*\omega_\pi)^{\ot n}
\ot \Oo(n\Delta)\ot d(\undE)^*\\
&=& \det E\boxt\det E^\vee(n\Delta) \ot \pi_2^* d(E)^*
\end{eqnarray*}
where $\pi_2$ is as in \eqref{pi-2}. This completes the proof of
the theorem.

\subsection{The determinant of the Szeg\"o kernel.} 
In the rank one case, \eqref{formula for Szego} expressed the Szeg\"o
kernel in terms of theta functions. In the higher rank case, however,
$\s_E$ is no longer a section of a line bundle, so cannot be expressed
directly this way. The Szeg\"o kernel $\s_E$ for
$E\in\Mf_X(n)\sm\Theta_X$ is a global section $$\s_E\in\Ho(\XX,E\boxt
E^{\vee}(\Delta))$$
(\defref{def of Szego}). On the other hand, by \propref{det prop} we
know that the pullback of the theta function by the difference map
$$\delta_E:\XX\to \Mf_X(n),\hskip.3in (x,y)\mapsto E\ot\Oo(y-x)$$ is a
section $$\delta_E^*\theta \in
\Ho(\XX,\det E\boxt\det\Ev(n\Delta))\ot\Theta_X|_E.$$ Thus it is
natural to relate $\delta^*_E\theta$ to a ``determinant'' of the
section $\s_E$.

It is convenient to interpret the section $\s_E$ of $p_1^*(E)\ot
p_2^*(\Ev)(\Delta)$, by dualizing the second factor, as a homomorphism
$$
\s_E\,\in\,
\Ho(\XX , \on{Hom}(p_2^*(E\ot \T_X)(-\Delta), p_1^* E))
$$
between two rank $n$ vector bundles on $\XX$. We may now take
its determinant as a homomorphism
$$
\on{det}(\s_E)\, \in\, 
\Ho(\XX ,\on{Hom}(p_2^*\on{det}(E\ot \T_X)(-n\Delta),
p_1^*\on{det}(E)))\, ,
$$
equivalently obtaining a section $\det\s_E\in
\Ho(\XX,\det E\boxt\det\Ev(n\Delta))$. We will find a cohomological
interpretation of this section, and of the nonabelian prime form,
which will make the comparison with the theta function immediate.

\subsubsection{\bf{Theorem}.}\label{kernel and theta}
 $\on{det}\s_E=\delta_E^*\theta/\theta(E)$ (under the identification
given by the nonabelian prime form).

\subsubsection{Proof.}
Recall that $\undE=p^*_X E$ denotes
the pullback of $E$ from $X$ to $X\times(\XX)$,
$\E=\undE(\Delta_{02}-\Delta_{01})$, and $\Pi:X\times(\XX)\to \XX$
denotes the projection.

First we rewrite the sheaves $p_2^*(E\ot \T_X)(-\Delta)$ and
$p_1^* E$ in terms of $\Pi$ and $\undE$. We have
\begin{eqnarray*}
p_1^*E&=& \Pi_* (\undE|_{\Delta_{01}})\\ p_2^*(E\ot
\T_X)(-\Delta)&=&\Pi_*
(\E/\undE(-\Delta_{01})).
\end{eqnarray*}
The second statement follows by observing that
$\undE(\Delta_{02}-\Delta_{01})/\undE(-\Delta_{01})$ is supported on
$\Delta_{02}$, and is given by
$$\undE(\Delta_{02})|_{\Delta_{02}}(-\Delta_{01}\cap \Delta_{02}),$$
which under the identification $\Pi|_{\Delta_{02}}:\Delta_{02}\cong
\XX$ and adjunction formula becomes $p_2^*(E\ot T_X)(-\Delta)$.

We next consider the long exact sequence of direct images, for the
projection $\Pi$, of the sequence of sheaves
\begin{equation}\label{first seq}
0\to \undE(-\Delta_{01})\to \E \to
\E/\undE(-\Delta_{01})\to 0\, .
\end{equation}
So we have
\begin{eqnarray*}
R^0\Pi_* \undE(-\Delta_{01}) \to &R^0\Pi_*\E& \to
R^0\Pi_* \E/\undE(-\Delta_{01})\to\\ 
R^1\Pi_* \undE(-\Delta_{01}) \to &
R^1\Pi_*\E& \to R^1\Pi_* \E/\undE(-\Delta_{01})
\end{eqnarray*}
The last term vanishes for dimension reasons, while the first term
vanishes since $E$ (and in particular $E(-x)$ for any $x\in X$) has no
global sections. Thus we obtain the following exact sequence:
\begin{equation}\label{theta sequence}
0\to R^0\Pi_*\E \to R^0\Pi_*
\E/\undE(-\Delta_{01})\stackrel{\pa_{\E}}\longrightarrow
R^1\Pi_* \undE(-\Delta_{01}) \to
R^1\Pi_*\E\to 0
\end{equation}

A similar consideration applies to the $\Pi$ pushforward of 
\begin{equation}\label{second seq}
0\to \undE(-\Delta_{01})\to \undE\to \undE|_{\Delta_{01}}\to 0,
\end{equation}
giving rise to
\begin{equation}\label{other theta sequence}
0\to R^0\Pi_*\undE\to R^0\Pi_*\undE|_{\Delta_{01}}
\stackrel{\pa_E}{\longrightarrow} R^1\Pi_* \undE(-\Delta_{01})\to
R^1\Pi_*\undE\to 0
\end{equation} In this case, however, both first and last terms
vanish so that $\pa_E$ is an isomorphism.

Thus we have a diagram 
\begin{equation}\label{cohomological Szego}
\begin{array}{ccc}
p_2^*(E\ot \T_X)(-\Delta)& \stackrel{\s_E}{\longrightarrow}& p_1^* E\\
\parallel&&\parallel\\
R^0\Pi_*\E/\undE(-\Delta_{01})&\stackrel{\pa_E\inv\circ
\pa_{\E}}{\longrightarrow}& R^0\Pi_* \undE|_{\Delta_{01}}
\end{array}
\end{equation}
which we claim commutes. It suffices to check that $\pa_E\inv \circ
\pa_\E$ restricted to the diagonal is the identity, since this
property characterizes $\s_E$ uniquely. However, over the diagonal
the two sequences \eqref{first seq} and \eqref{second seq} coincide,
so that $\pa_\E|_\Delta=\pa_E|_\Delta$.

The theorem now follows from the cohomological description of theta
functions. 
Namely, the complexes
$$R^0\Pi_*
\E/\undE(-\Delta_{01})\stackrel{\pa_{\E}}\longrightarrow
R^1\Pi_* \undE(-\Delta_{01})$$
and
$$R^0\Pi_*\undE|_{\Delta_{01}} \stackrel{\pa_E}{\longrightarrow}
R^1\Pi_* \undE(-\Delta_{01})$$ (from \eqref{theta sequence},
\eqref{other theta sequence}) are easily seen to represent the
determinant bundle for the coherent sheaves $\E$ and $\undE$ on
$X\times(\XX)$, in the sense of \secref{det bundles}. For example,
for the first complex $\K^\cdot$ we must check that for any coherent
sheaf $\mc G$ on $\XX$ we have
$$R^i\Pi_*(\E\ot\Pi^*{\mc G})\cong {\rm H}^i(\K^{\cdot}\ot {\mc G}).$$
Thus we consider as before (when $\mc G$ was $\Oo_{\XX}$) the long exact
sequence associated to \eqref{first seq} tensored with $\Pi^*{\mc G}$.
But by the projection formula the first and last terms of the
resulting six--term sequence will again vanish, and the two middle
terms will give our $\K^0\ot{\mc G}$ and $\K^1\ot{\mc G}$, as desired.

It follows that we have a canonical isomorphism of line bundles
$$\un{\on{Hom}}_{\XX}(p_2^*\on{det}(E\ot \T_X)(-n\Delta),
p_1^*\on{det}(E))\cong d(\E)^*\ot d(\undE),$$ (where $d(\undE)$ is the
trivial line bundle $d(E)\ot_{\C}\Oo_{\XX}$) and that
$\det\pa_\E=\delta_E^*\theta$ and $\det\pa_E$ is the
constant $\theta(E)$.
Such an isomorphism was
constructed in \thmref{det theorem}, the nonabelian prime form. Since
$\XX$ is projective, the isomorphisms must be proportional. In fact,
the isomorphisms are equal, since along the diagonal both restrict to
the identity map $\Oo\to d(E)^*\ot d(E)$. This gives the cohomological
description of the prime form.
Tracing through these identifications we obtain
$$\on{det}\s_E=\delta_E^*\theta/\theta(E),$$ as desired.

\section{Connections and the Szeg\"o kernel}\label{conns and Szego}

\subsection{Nonabelian Szeg\"o Kernels.}
For our purpose it is convenient to introduce the subsheaf
$\No(E)\subset \M(E)$
consisting of all sections whose restriction to the diagonal
$\Delta$ is a scalar multiple of the identity,
$$\No(E)\,=\,\{s\in \M(E) \big\vert\,\, 
s|_{\Delta}\in \Oo_{\Delta}\cdot
\on{Id}\}.
$$
In other words, $\No(E)$ is the inverse image of the line-subbundle of
$\text{End}\,E=\M(E)|_\Delta$ defined by scalar operators. For $E$ of
rank $1$ we have $\No(E)=\M(E)$.

\subsubsection{\bf{Proposition}.}\label{kernels}
\begin{enumerate}
\item If ${\rm h}^0(E)={\rm h}^1(E)=0$, then $\Ho(X\times
X,\No(E))=\C\cdot\s_E$, where $\s_E$, the Szeg\"o kernel of $E$,
is the unique section with $\s_E|_{\Delta}=\on{Id}$.
\item Otherwise, the inclusion $\Ho(\XX,E\boxt
E^{\vee})\hookrightarrow \Ho(\XX,\No(E))$ is an isomorphism,
i.e., all global sections of $\No(E)$ vanish on $\Delta$. 
\end{enumerate}

\subsubsection{Proof.}
Consider the exact sequence of cohomology
\begin{equation*}
0\to \Ho(\XX, E\boxt \Ev)\to
\Ho(\XX, \No(E))\to
\Ho(\Delta,\Oo_\Delta)\stackrel{H}{\to}
\Hi(\XX, E\boxt \Ev)
\end{equation*}
obtained from the short exact sequence
\begin{equation}\label{sequence for No}
0\, \longrightarrow \, E\boxt \Ev \, \longrightarrow \,
\No(E) \, \longrightarrow \, {\mathcal O}_\Delta \, \longrightarrow
\, 0
\end{equation}
In the projection $\No(E)\vert_\Delta \longrightarrow
{\mathcal O}_\Delta$, the section of ${\mathcal O}_\Delta$
defined by the constant
function $1$ corresponds to the identity automorphism of $E$.

By the K\"unneth formula, the cohomology 
$$\Hi(\XX,E\boxt\Ev)=\left(\Ho(X,E)\ot \Hi(X,\Ev)\right)\oplus
\left(\Hi(X,E)\ot \Ho(X,\Ev)\right).$$ This space has a canonical
element coming from the Serre duality pairing between the cohomologies
of $E$ and $\Ev$. It is straight-forward to check that the image
$H(1)$ of the function $1\in \Ho(\Oo_{\Delta})=\C$
(in the above exact sequence of cohomology) is this
Serre duality element. It follows that the homomorphism
$H$ is an embedding whenever
either $\Ho(X,E)$ or $\Hi(X,E)$ is nonzero, while if both vanish then
$H=0$. The proposition follows immediately.

\subsubsection{Remark.}\label{Serre}
This construction is a special case of the description of Serre
duality via residues at the diagonal for any smooth
projective variety $X$ of dimension $d$. Consider the map
$$
{\rm H}^d_{\Delta}(\XX, E\boxt \Ev)\longrightarrow {\rm H}^d(\XX,
E\boxt \Ev)=\bigoplus_{0\leq i\leq d}
{\rm H}^i(X,E)\ot {\rm H}^{d-i}(X,E^{\vee})
$$
{}from top local cohomology at the diagonal to middle cohomology on
$\XX$. The local cohomology is identified by the Grothendieck--Sato
formula with the differential operators $\Ho(X,\D(E,E))$ from $E$ to
$E$. The image of the identity operator is the Serre duality pairing.

\subsubsection{Remark: Even nonsingular theta characteristics.}
If the theta characteristic $\Ohalf$ lies in the complement
$\Pg(X)\sm\Theta_X$, in other words $\Ohalf$ is even and nonsingular,
then the Szeg\"o kernel $\s$ associated to $\Ohalf$ is the classical
Szeg\"o kernel introduced in \cite{Hawley}. Since the transpose of the
Szeg\"o kernel for $E$ clearly agrees with the Szeg\"o kernel for
$E^{\vee}$, it follows that $\s$ is symmetric. Hence its restriction
to $2\Delta$ necessarily agrees with the canonical trivialization
$\mu_1$ of $\M_1|_{2\Delta}=\Ohalf\boxt\Ohalf(\Delta)|_{2\Delta}$.

\subsection{An identity for the Szeg\"o kernel}\label{identity}
Let $E$ be a vector bundle over $X$ with $h^0(E)=h^1(E)=0$, and
Szeg\"o kernel $\s_E$. Let $\s_E^{ij}$ denote the pullback of $\s_E$
to $X\times\XX$ along the $i$th and $j$th factors. Regarding $\s_E$ as
a homomorphism from $p^*_2E$ to $E\boxt\Omega_X(\Delta)$
over $X\times X$, the
composition $\s_E^{01}\circ \s_E^{12}$ gives a vector bundle
homomorphism from the pullback $p^*_3E$ along the third factor to
$E\boxt\Omega_X\boxt\Omega_X(\Delta_{01}+\Delta_{12})$; here
$p_3\, :\, X\times X\times X\, \longrightarrow\, X$
be the projection to the third factor.
We will consider $\s_E^{01}\circ \s_E^{12}$ as a meromorphic
homomorphism from $p_3^*E$ to $E\boxt \Omega_X\boxt\Omega_X$, with
pole of order one along the 
(connected reduced) divisor $\Delta_{01}+\Delta_{12}$. 

Let $f\, :\, X\, \longrightarrow\, {\mathbb C}\cup\{\infty\} \,=
\,{\mathbb C}{\mathbb P}^1$ be a nonconstant meromorphic
function. Take a point $c\in {\mathbb C}$ such that $f^{-1}(c)$ is
reduced. In other words, if the degree of $f$ is $l$, then $f^{-1}(c)$
consists of $l$ distinct points. Without loss of generality
(translating $f$ if necessary) we assume $c=0$.  The following
proposition describes an identity satisfied by the Szeg\"o kernel.

\subsubsection{\bf{Proposition}.}\label{Szego-identity}
For any $(x,y)\in X\times X \setminus \Delta$, the identity
$$
\sum_{\alpha \in f^{-1}(0)}\frac{1}{df (\alpha)}
\s_E (x,\alpha)\circ \s_E (\alpha,y) \, =\,
\frac{f(y)-f(x)}{f(y)f(x)}\s_E (x,y)
$$
is satisfied.

\subsubsection{Proof.}
Fix a point $y\in X$ and fix a vector $e\in E_y$ in the fiber.
Also, fix a nonzero holomorphic tangent vector $v\in T_yX\setminus
\{0\}$. Now, the map defined by
$$
x \, \longmapsto \,
\langle \sum_{\alpha \in f^{-1}(0)}\frac{1}{df (\alpha)}
\s_E (x,\alpha)\circ \s_E (\alpha,y)(e)\, ,v\rangle
$$
is a meromorphic section of $E$. Note that the condition on
$0$ ensures that $df (\alpha)\not= 0$. Therefore,
$df (\alpha)$ and $v$ trivialize the fibers $(\Omega_X)_\alpha$
and $(\Omega_X)_y$ respectively; here
$\langle-\, ,v\rangle$ denotes the contraction
of $(\Omega_X)_y$ by $v$. We will denote this
meromorphic section of $E$ by $A(e,v)$.

Similarly,
$$
x \, \longmapsto \,
\langle \frac{f(y)-f(x)}{f(y)f(x)}\s_E (x,y)\, ,v\rangle
$$
is a meromorphic section of $E$, which we will denote by
$B(e,v)$.

The proposition is equivalent to the assertion
$A(e,v)\, =\, B(e,v)$ for all $e,v$.
Since $E$ does not admit any nonzero
holomorphic section, it suffices to show that the poles of
$A(e,v)$ and $B(e,v)$ coincide. Indeed, in that case $A(e,v)-
B(e,v)$ is a holomorphic section and thus must be identically zero.

First consider the case $y\,\notin\, f^{-1}(0)$. Then the poles
of the section $A(e,v)$ are at $f^{-1}(0)$, each pole is
of order one and the residue at $\alpha \in f^{-1}(0)$ is

$$
\frac{\langle \s_E(\alpha,y)\, ,v\rangle}{df(\alpha)}
\, \in\, E_x\otimes T_\alpha X\, .
$$
Since ${\mathcal O}_X(\alpha)$
is identified (by the adjunction formula) with $T_\alpha X$,
the residue is an element of $E_x\otimes T_\alpha X$;
here $1/df(\alpha)$ denotes the element of $T_\alpha X$
dual to the nonzero element $df(\alpha)$ of $(\Omega_X)_\alpha$.
That the residue of $A(e,v)$ coincides with $\langle \s_E(\alpha,y)\,
,v\rangle/df(\alpha)$ is an immediate consequence of the property
of $\s_E$ that its restriction to $\Delta$ is the identity
automorphism of $E$.

The poles of $B(e,v)$ are contained in
$f^{-1}(0)\cup\{y\}$. Now, $y$ is a removable singularity
since the pole of $\s_E$ at $y$ is canceled by the vanishing
of $f(y)-f(x)$. Clearly, the pole of $B(e,v)$ at any $\alpha
\in f^{-1}(0)$ is of order one, and the residue 
is $\langle \s_E(\alpha,y)\, ,v\rangle/df(\alpha)$.
(First note that the numerator $f(y)-f(x)$ cancels with
$f(y)-0$ and then the expression of residue follows from the fact
that the residue of $1/f(x)$ is
$1/df(\alpha)$.)

Since the residues of $A(e,v)$ and $B(e,v)$ coincide,
the proposition is proved under the assumption that
$y\,\notin\, f^{-1}(0)$. If $y\,\in\, f^{-1}(0)$, then
it follows by continuity as $y$ is in the closure
of the complement of $f^{-1}(0)$. (This case can also be proved
by exactly the same way as done above for
$y\,\notin\, f^{-1}(0)$. The
only point to take into account is that using the tangent
vector $v$ the line $({\mathcal O}_{X\times X}(\Delta))_{(y,y)}$
has to be trivialized.)

\subsubsection{Remark: The degenerate case of the
Szeg\"o kernel identity.}
\propref{Szego-identity} can be extended to the case
$y=x$, that is, across the diagonal. Taking the limit
as $x$ approaches $y$, the identity in \propref{Szego-identity}
becomes
$$
\sum_{\alpha \in f^{-1}(0)}\frac{1}{df (\alpha)}
\s_E (y,\alpha)\circ \s_E (\alpha,y) \, =\,
\frac{df(y)}{f(y)^2}\text{Id}_{E_y}\, .
$$
Indeed, this is an immediate consequence of the fact that
the evaluation
$\s_E (y,y)$ is the identity automorphism of $E_y$
(see Definition \ref{def of Szego}).

\subsection{The sheaf of Szeg\"o kernels.}
We would like to describe a version of \propref{kernels} with
parameters, in other words to let the bundle $E$ (and possibly the
curve $X$) vary. For this purpose we require a universal bundle to
substitute for $E\boxt\Ev$. Such a universal bundle is available
(tautologically) over the moduli {\em stack} of bundles, so all of
our arguments will work in the stack setting. Alternatively, for
concreteness we restrict to the stable locus $\Mf_X^s(n)\subset
\Mf_X(n)$, recalling that the codimension of the strictly semistable
locus is greater than one, so that line bundles and their sections are
determined by the restriction to the stable locus. Recall
(\secref{moduli spaces}) that there is a universal bundle
$\Ef\boxt\Ef^{\vee}$ over $X\times X\times\Mf_X^s(n)$ whose
restriction to $X\times X\times\{E\}$ is $E\boxtimes E^*$. Using this
bundle, we define a canonical vector bundle $\No\subset
\Ef\boxt\Ef^{\vee}(\Delta)$ over $X\times X\times\Mf_X^s(n)$ whose
restriction to any $X\times X\times\{E\}$ is $\No(E)$.

Let $\pi_2\, :\, \XX\times\Mf_X^s(n)\,\longrightarrow\,
 \Mf_X^s(n)$ be the natural projection (as in \eqref{pi-2}).
We now describe the direct image sheaf
$$
{\mf S}\,=\,\pi_{2*}\No
$$
on $\Mf_X(n)$, where the Szeg\"o kernel naturally lives,
for the above projection $\pi_2$.

\subsubsection{\bf{Proposition}.} \label{normalized Szego}
$\mf S$ is naturally a torsion--free subsheaf of the determinant line
bundle $\Theta_X^*$. The inclusion is an isomorphism on
$\Mf_X^s(n)\sm\Thetasing$, where its inverse is given by the {\em
normalized Szeg\"o kernel} $\sbar:=\s\cdot\theta$, a nowhere vanishing
section of ${\mf S}\ot\Theta_X$ off $\Thetasing$.

\subsubsection{Proof.}
Consider the relative version
$$0\to \Ef\boxt\Ef^{\vee}\to \No\to\Oo_\Delta\to 0$$
of \eqref{sequence for No}, and push it forward
along $\pi_2$, obtaining the
exact sequence 
$$
0\longrightarrow \pi_{2*}(\Ef\boxt\Ef^{\vee})\longrightarrow {\mf
S}\stackrel{\rm{Res}}{\longrightarrow} \Oo\longrightarrow\cdots
$$ 

The sheaf $\pi_{2*}(\Ef\boxt\Ef^{\vee})$
vanishes off of $\Theta_X$, and is a line
bundle on $\Theta_X^{\circ}$.
In fact, for $E\in\Theta_X^{\circ}$, 
the vector spaces $\Ho(X,E)$ and $\Hi(X,E)^*\cong \Ho(X,E^\vee)$ are
one--dimensional, so that we have isomorphisms
$$\Ho(\XX,E\boxt\Ev)\cong\Ho(X,E)\ot\Ho(X,E^{\vee}) \cong
\Theta_X^*|_E.$$ As $E$ varies over $\Theta_X^\circ$ we obtain an
isomorphism $\Theta_X^*|_{\Theta_X^\circ}\cong
\pi_{2*}(\Ef\boxt\Ef^{\vee})$.

Since $\No$ is a locally free sheaf, ${\mf S}\,=\,\pi_{2*} \No$ is
torsion--free. On $\Mf_X^s(n)\sm\Theta_X$, the morphism $\on{Res}$
above is an isomorphism. Moreover, $\mf S$ remains locally free on
$\Mf_X(n)\sm\Thetasing$. To see this observe (by \propref{kernels}
and the above discussion) that on $\Theta_X^{\circ}$, $\mf S$ is
identified with the line bundle $\Theta_X^*$, so (by Nakayama's lemma)
it is locally free of rank one near any $E\in\Theta_X^{\circ}$.

The morphism $\on{Res}$ is a regular section of the line bundle
$\un{\on{Hom}}({\mf S},\Oo)$. It is nonvanishing off $\Theta_X$, and
vanishes to first order along $\Theta_X$. This identifies
$\un{\on{Hom}}({\mf S},\Oo)$ with $\Theta_X$ and provides the desired
embedding ${\mf S}\subset \Theta_X^*$, an isomorphism off of
$\Thetasing$. Moreover, since by definition the residue of $\s_E$ is
the identity for every $E$, tracing through the isomorphism we see
that the section $\sbar=\s\ot\theta$ of ${\mf S}\ot\Theta_X=
\un{\on{Hom}}(\Theta_X^*,{\mf S})$ on $\Mf_X^s(n)\sm\Theta_X$
is identified with the inverse of the identification above. In
particular this section extends across $\Theta_X^{\circ}$ in a
nonzero fashion.

\subsubsection{Remark.} Note that the statement 
(and proof) above generalize immediately when the curve $X$
and bundle are allowed to vary over the moduli space $\Mf_g(n)$.

\subsection{Connections on $\Theta_X$.}
Consider the affine space $\Conn_X(E_0)$ of twisted connections on $E$.
This is a torsor for the space of Higgs fields, which is the cotangent
space to $\Mf_X(n)$ at $E$ (at smooth points -- in particular for $E$
stable). Thus as $E$ varies, the sheaf on $\Mf_X(n)$ defined by
$\Conn_X(E_0)$ forms an $\Omega$--torsor. This $\Omega$--torsor will be
denoted by $\Conn_X(n)$. On the complement $\Mf_X(n)\sm\Theta_X$, the
restriction of the Szeg\"o kernel $\s_E|_{2\Delta}$ to $2\Delta$ gives
a canonical section of $\Conn_X(n)$, denoted $\s|_{2\Delta}$.

\subsubsection{\bf{Theorem}.}\label{dlog theta}
There is a unique isomorphism of $\Omega$--torsors
$\Conn_X(n)\cong \Conn_X(\Theta)$ on $\Mf_X(n)$ sending
$\s|_{2\Delta}$ to $\dlog \theta$.

\subsubsection{Proof.}
The Szeg\"o kernel $\s|_{2\Delta}$ is a meromorphic section of the
$\Omega$--torsor $\Conn_X(n)$, trivializing it on
$\Mf_X(n)\sm\Theta_X$. Moreover it has only a first order pole along
$\Theta_X$, in fact as we will see a {\em logarithmic} pole. It follows
therefore that the class of $\Conn_X(n)$ in $\Hi(\Mf_X(n),\Omega^1)$ is
a complex multiple of the divisor class $[\Theta_X]$, with the
multiplicity given by the residue of $\s|_{2\Delta}$ along $\Theta_X$.
(Note that it makes sense to speak of the residue of a section of
any $\Omega$--torsor with first--order poles, since the expressions
of this section in any two local trivializations differ by a regular
one--form.) Hence we need to check $\s|_{2\Delta}$ is logarithmic,
calculate its residue and compare it to that of $\dlog\theta$.

We first describe the conormal bundle $N^*_{\Theta_X^{\circ}}$ to
$\Theta_X^{\circ}$. By the adjunction formula, the conormal bundle is
isomorphic to the restriction of the dual line bundle $\Theta_X^*$. This
isomorphism is given by $\dlog\theta|_{\Theta_X^\circ}$, a nonvanishing
section of $N^*_{\Theta_X^\circ}\ot\Oo(\Theta_X)|_{\Theta_X^\circ}
\subset\Omega^1_{\Mf_X(n)}(\Theta_X)|_{\Theta_X^{\circ}}$. (This is
also identified with the one--jet $J^1\theta$, which along the divisor
of $\theta$ is a section of $\Omega^1_{\Mf_X(n)}\ot\Oo(\Theta_X)$.)

On the other hand, identifying $\Theta_X^{\circ}$ with the
Brill--Noether moduli space of bundles with a single holomorphic
section, standard deformation theory identifies its  tangent space
at $E\in\Theta_X^{\circ}$ with the kernel of the cup product map
\begin{equation*}
\Hi(X,\on{End}E)\longrightarrow \on{Hom}(\Ho(X,E),\Hi(X,E)).
\end{equation*}
(Given a section of $E$ and a first order deformation of the bundle,
we obtain a class in $\Hi(X,E)$, while a compatible deformation of the
section makes this class into a coboundary.)  Dually, the conormal
bundle is identified with the image of the Petri map
\begin{equation*}
\Ho(X,E)\ot\Ho(X,\Ev)\, \longrightarrow\, \Ho(X,E\ot\Ev)\, ,
\end{equation*}
or equivalently the restriction to the diagonal
$$\Ho(\XX,E\boxt\Ev)\to\Ho(X,\on{End}E\otimes \Omega_X)\, ,$$ which are
injective when ${\rm h}^0(E)={\rm h}^0(\Ev)=1$. These two
identifications of the
conormal bundle are compatible with the natural isomorphism
$\Theta_X^*|_E=\Ho(X,E)\ot\Ho(X,\Ev)$ for $E\in\Theta_X^\circ$.

For $E\in\Theta_X^\circ$ the restriction of the normalized Szeg\"o
kernel $\sbar|_{2\Delta}$ is a section of
$E\boxt\Ev(\Delta)|_{2\Delta}\ot\Theta_X\vert_{\{E\}}$ vanishing on
the diagonal. In other words, it is a $\Theta_X|_{\{E\}}$--twisted
section of $E\boxt\Ev|_\Delta=\on{End}E\ot\Omega$. Moreover, this
twisted Higgs field is the image of
$$
\sbar_E\in \Ho(\XX,E\boxt\Ev)\ot\Theta_X|_{\{E\}}\, ,
$$
hence (by the above
identifications) it is a nonvanishing function $\sbar|_{2\Delta}\in
N^*_{\Theta_X^{\circ}}\ot\Theta_X\cong \Oo$ along $\Theta_X^\circ$.

It follows that $\s|_{2\Delta}$ has log--poles along $\Theta_X^\circ$,
since its singular part is contained in the conormal direction.
Therefore its residue can
be read of as the function on $\Theta_X^\circ$
given by the restriction $\sbar|_{2\Delta}$. But by
\propref{normalized Szego}, the normalized Szeg\"o kernel (along
$\Theta_X^\circ$) is simply the expression of the canonical isomorphism
$\Theta_X^*|_E=\Ho(X,E)\ot\Ho(X,\Ev)$. In other words, the residue of
$\s|_{2\Delta}$ along $\Theta_X$ is $1$, and the class of the
$\Omega$--torsor $\Conn_X(n)$ is the (Chern) class of $\Theta_X$.

Consider the meromorphic identification of $\Omega$--torsors $\Conn_X(n)
\to\Conn_X(\Theta)$ defined by $\s|_{2\Delta}\mapsto \dlog\theta$.
This isomorphism is clearly regular outside $\Theta_X$, and the above
identification of the
singular parts shows that it remains regular along
$\Theta_X$. Thus we have constructed a global and uniquely characterized
isomorphism of $\Omega$--torsors, as desired.

\subsubsection{Uniqueness and symplectic structure.} \label{uniqueness}
As explained in \secref{rigidity}, $\Omega$--torsors on moduli spaces
enjoy a remarkable rigidity. Since the isomorphism in \thmref{dlog
theta} is valid over arbitrary families of curves, and there are no
one--forms on $\Mf_g(n)$ relative to $\Mf_g$, it follows that this
isomorphism agrees with the isomorphisms described in \cite{Fa} and
\cite{BS}. The latter isomorphisms are in fact isomorphisms of {\em
twisted cotangent bundles}, i.e., they carry the symplectic structure
on flat connections (given by the cup product on de Rham cohomology of
a connection) to that on $\Conn_X(\Theta)$ (see
\secref{symplectic}). While we have not addressed the relation of
Szeg\"o kernels to the natural symplectic structure on $\Conn_X(n)$
above, this relation follows immediately by comparison with these
other constructions:

\subsubsection{\bf{Corollary}.}\label{Lagrangian} The section
of $\Conn_X(n)$ over 
$\Mf_X(n)\sm\Theta_X$ given by the Szeg\"o kernel is Lagrangian with
respect to the natural twisted cotangent bundle symplectic structure
on the moduli space of connections.

\section{Extended Connections}\label{exconns and Szego}

\subsection{Extended Higgs Fields.}
Let $\Mf_g(n)$ denote the moduli space of semistable bundles of rank
$n$ and Euler characteristic zero over the moduli space of curves
$\Mf_g$. This is an orbifold (Deligne--Mumford stack) which fibers
over the moduli space of curves, with fiber at $X$ the moduli space
$\Mf_X(n)$. It carries a universal theta divisor, also denoted by
$\Theta_X$, which specializes to the theta divisor on each $\Mf_X(n)$.
(See \cite{universal} for a discussion of this moduli space.)

The tangent space to $\Mf_g(n)$ at a smooth point $(X,E)$ is
calculated, by standard deformation theory, to be
${\rm H}^1(X,\on{At}(E))$,
where $\on{At}(E)$ is the Atiyah bundle of $E$. The cotangent space
to $\Mf_g(n)$ may thus be described concretely as an $\Ho$ of the dual
sheaf of kernels. Using the residue pairing along the diagonal and the
trace on $\on{End}E$, the result is (after transposition of factors)
the following (see \cite{BS}):

\subsubsection{Definition.}
The space $\exhiggs_X(E)$ of extended Higgs fields on $E$ is
the sheaf 
$$
\exhiggs_X(E)\,=\,\{\phi\in E\boxt
E^*(-\Delta)|_{2\Delta}\}/(\on{Ker}(\on{tr})\subset E\boxt
E^*(-2\Delta)|_{\Delta})
$$
supported over $2\Delta$, where
$$
\on{tr}\, :\, E\boxt E^*(-2\Delta)\vert_\Delta\, =\,
\on{End} E\otimes \Omega^{\otimes 2}_X\, \longrightarrow\,
\Omega^{\otimes 2}_X
$$
is the trace map. (Note that $E\boxt E^*(-2\Delta)\vert\Delta\,
\subset\, E\boxt E^*(-\Delta)|_{2\Delta}$ is the subsheaf
vanishing on $\Delta\, \subset\, 2\Delta$.)

The direct image $p_{1*}\exhiggs_X(E)$ is a vector bundle of
rank $n^2+1$ over $X$. This vector bundle will also be called
the extended Higgs fields on $E$ (as we have a natural
identification defined by taking this direct image).

\subsubsection{} It is easy to check that the bundle
$\exhiggs_X(E)$ is an extension
\begin{equation}\label{exhiggs equation}
0\to \Omega^{\otimes 2}_X\to \exhiggs_X(E)\to \on{End}E\ot \Omega_X\to 0
\end{equation}
of the sheaf of Higgs fields $\on{End}E\ot \Omega_X$ by the
sheaf of quadratic differentials $\Omega^{\otimes 2}_X$.

\subsubsection{\bf{Proposition}.} The cotangent
space to $\Mf_g(n)$ at any point $(X,E)$ is
given by the space of global extended Higgs fields on
$E$, namely $\Ho(X,\exhiggs_X(E))$.

\subsubsection{Proof.}
The tangent space to $\Mf_g(n)$ at any point $(X,E)$
is given by $\Hi(X,\text{At}(E))$, where $\text{At}(E)$
is the Atiyah bundle defined in \secref{Atiya-bundl.-def.}
(see \cite{BS}). The Serre duality gives
$\Hi(X,\text{At}(E))^*\, \cong\, \Ho(X,\text{At}(E)^*\otimes
\Omega_X)$. The proof of the proposition will be completed
by giving a canonical isomorphism of
$\text{At}(E)^*\otimes\Omega_X$ with the vector bundle
$\exhiggs_X(E)$.

To construct the isomorphism, first note that the pulled back
line bundle $(p^*_1\Omega_X)\vert_{2\Delta}$ over
$2\Delta$ is naturally identified with the line bundle
$\Oo_{X\times X}(-\Delta)\vert_{2\Delta}$. To construct this
isomorphism, observe that a holomorphic function $f$
defined on a connected
symmetric analytic open subset $U\, \subset\,
X\times X$ (by symmetric we mean invariant under the
involution $\sigma$ of $X\times X$ that interchanges factors)
satisfying the identity $f\circ \sigma \, =\, -f$ has the
following property: if $f$ vanishes on $\Delta\cap U$
of order at least two, then $f$ must vanish
on $\Delta\cap U$ of order at least
three. This observation gives an isomorphism of
$\Oo_{X\times X}(-\Delta)\vert_{2\Delta}$ with
$(p^*_1\Omega_X)\vert_{2\Delta}$. Indeed, taking any function
$f$ on such a symmetric open subset $U$ that vanishes
on $\Delta\cap U$ exactly of order one, project the
section ${\rm d} f$ of ${\Omega}^1_U$ to a
section of $(p^*_1\Omega_X)\vert_U$, where the projection is
defined using decomposition ${\Omega}^1_{X\times X}\,
\cong\, p^*_1\Omega_X\oplus p^*_2\Omega_X$. The resulting
isomorphism of $\Oo_{X\times X}(-\Delta)\vert_{2\Delta}$ with
$(p^*_1\Omega_X)\vert_{2\Delta}$ does not depend on the
choice of $f$. This is an immediate consequence of the
above observation.

{}From the definition of $\text{At}(E)$
in \secref{Atiya-bundl.-def.} it follows immediately that
$\text{At}(E)^*$ is a quotient of $J^1(E)\otimes E^*$. Recall
that $J^1(E)$ is, by definition, the direct
image $p_{2*} ((p^*_1 E)\vert_{2\Delta})$. Consider the
quotient sheaf
$$
{\mathcal V} \,:=\,\{\phi\in E\boxt
E^*\otimes p^*_1\Omega_X|_{2\Delta}\}/(\on{Ker}(\on{tr})\subset
E\boxt E^*\otimes p^*_1\Omega_X(-\Delta)|_{\Delta})
$$
on $2\Delta$, where
$$
\on{tr}\, :\, E\boxt E^*\otimes p^*_1
\Omega_X(-\Delta)\vert_\Delta\,
=\, \on{End} E\otimes \Omega^{\otimes 2}_X\, \longrightarrow\,
\Omega^{\otimes 2}_X
$$
is the trace map. Using the above remarks it follows that
$p_{1*}{\mathcal V} \, \cong\, \text{At}(E)^*\otimes\Omega_X$.

Finally, using the earlier obtained isomorphism of
$\Oo_{X\times X}(-\Delta)\vert_{2\Delta}$ with
$(p^*_1\Omega_X)\vert_{2\Delta}$, the vector bundle
$\text{At}(E)^*\otimes\Omega_X$ over $X$ is identified
with $\exhiggs_X(E)$. This completes the proof of the
proposition.

\subsection{Extended connections.}
We introduce the space of extended connections, a torsor over
the space of extended Higgs fields, using the language
of kernels. 

Recall that a connection on a vector bundle $E$ may be realized,
for any $\nu $, by a section of $\M_\nu(E)$ on $2\Delta$, whose
restriction to $\Delta$ is the identity endomorphism of $E$. It is
convenient to describe the sheaf of sections on $2\Delta$ as the
quotient of $\M_\nu(E)$ by sections that vanish on $2\Delta$: we have
an identification
$$
\Conn_X(E)\cong\{\varphi\,\in\, \M_\nu(E)/\M_\nu(E)(-2\Delta)\,
\big\vert\, \varphi|_\Delta =\on{Id}_E \}.
$$

In order to define extended connections, we would like to extend the
connection kernels on $2\Delta$ by {\em scalar} operators on
$3\Delta$. Thus instead of
quotienting out $\M_\nu(E)$ by {\em all} sections
that vanish on $2\Delta$, we quotient out by those whose leading term
is traceless. 
Let $\on{End}^0(E)\, \subset\, \on{End}E$ denote the
subbundle of traceless endomorphisms of $E$ and define
$$
\M_\nu(E)(-2\Delta)^0 \, :=\, \Omega_X^{\ot 2}\otimes\on{End}^0(E)
\, \subset\, \M_\nu(E)(-2\Delta)|_{\Delta} \,=\,
\Omega_X^{\otimes 2}\otimes \on{End}E \, .
$$
So $\M_\nu(E)(-2\Delta)^0$ coincides with the kernel of the
homomorphism
$$
\M_\nu(E)(-2\Delta)|_{\Delta}\, =\,
\Omega_X^{\otimes 2}\otimes \on{End}E\,
\stackrel{\on{tr}_E}{\longrightarrow}\, \Omega_X^{\ot 2}.
$$
(Note that $\M_\nu(E)(-2\Delta)|_{\Delta}$ is identified with the
sections of $\M_{\nu}(E)|_{3\Delta}$ vanishing on $2\Delta$.)

\subsubsection{Definition.}\label{ex.con.def.}
The sheaf of extended connections is
$$
\excon_X^\nu(E)=\{\vphi\,\in\, 
\M_\nu(E)/\M_\nu(E)(-2\Delta)^0\,\big\vert
\, \vphi|_\Delta =\on{Id}_E\},
$$
the sections of the quotient of
$\M_\nu(E)$ over $3\Delta$ by traceless sections
$\M_\nu(E)(-2\Delta)^0$, whose restriction to $\Delta$ is the
identity automorphism of $E$.

\subsubsection{}
It follows that restriction to $2\Delta$ defines a map $\Pi_\nu
:\excon_X^\nu(E)\to \Conn_X(E)$ to the space of connections on $E$.
If we take $\vphi_1,\vphi_2\,\in\, \excon_X^\nu(E)$
defining the same connection on $E$, that is
$\vphi_1|_{2\Delta}=\vphi_2|_{2\Delta}$, then their difference $\vphi_1
-\vphi_2=q\on{Id}_E$ with $q$ a quadratic differential on $X$. Thus
restriction to $2\Delta$ makes $\excon_X^\nu(E)$ an affine bundle for
$\Ho(X,\, \Omega_X^{\ot 2})$ over $\Conn_X(E)$.
 
\subsubsection{\bf{Proposition}.}\label{torsor for excons}
For every $\nu\in\Z$, the space $\excon_X^\nu(E)$ is naturally an
affine space for $\exhiggs_X(E)$. 

\subsubsection{Proof.} 
We claim that for every $\nu$, the vector bundle
\begin{equation*}
\exhiggs_X^\nu(E)=\{\vphi\,\in\, \M_\nu(E)/\M_\nu(E)(-2\Delta)^0\,
\big\vert \, \vphi|_\Delta=0\}
\end{equation*}
is isomorphic to $\exhiggs_X(E)$.
The isomorphism is given by tensoring with the section
$\mu_{-\nu}$, i.e., by the identification
$$\Gamma(\M_\nu(-\Delta)|_{3\Delta})=\Gamma(\Omega_X^{\nu/2}\boxt
\Omega_X^{\nu/2}((i-1)\Delta)|_{2\Delta})
\stackrel{\ot\mu_{-\nu}}{\longrightarrow}
\Gamma(\Oo(-\Delta)|_{2\Delta}).$$ It is clear that $\excon_X^{\nu}(E)$
is an affine bundle for the vector bundle
$\exhiggs_X^\nu(E)$, and hence for $\exhiggs_X(E)$.

\subsection{Varying the curve.}

We would like to use the Szeg\"o kernel $\s_E$ to define a section of
a twisted cotangent bundle over the moduli space $\Mf_g(n)$, extending
the construction $\s|_{2\Delta}$ along $\Mf_X(n)$. Thus we introduce
the extended analogue of $\Conn_X(E_0)$: the space $\excon_X(E_0)$ of
twisted extended connections on $E$ is
$$\excon_X(E_0)=\{s\in \Ho(3\Delta,\M(E)|_{3\Delta}) \,\vert\,
s|_\Delta=\on{Id}_E\}/ (\on{Ker}(\on{tr})\subset E\boxt
E^{\vee}(-\Delta)|_{3\Delta}).$$ Again this space depends only on $E$
and not on the choice of theta characteristic $\Ohalf$.

Now consider the projection of $\s_E|_{3\Delta}$ to $\excon_X(E_0)$,
which we also denote by $\s_E|_{3\Delta}$. This defines a section of
the $\Omega$--torsor $\excon_g(n)$ over $(X,E)\in\Mf_g(n)\sm\Theta_g$.

\subsubsection{\bf{Theorem}.}\label{dlog theta over curves}
There is a unique isomorphism of $\Omega$--torsors
\begin{equation*}
\excon_g(n)\, \cong\, \Conn_g(\Theta_g)
\end{equation*}
over $\Mf_g(n)$ sending $\s|_{3\Delta}$ to $\dlog \theta$.

\subsubsection{Proof.} 
The proof is a direct generalization of that of \thmref{dlog theta}.
Namely, we claim that $\s|_{3\Delta}$ is a section of $\excon_g(n)$
with logarithmic poles along the universal theta divisor $\Theta_g$,
and with residue $1$. As was the case for fixed curve, the logarithmic
derivative of the theta function identifies the conormal bundle to
$\Theta_g^\circ$ with the restriction of $\Theta_g^*$.  Deformation
theory now identifies the tangent space to $\Theta_g^\circ$ with the
kernel of a natural map
\begin{equation*}
\Hi(X,\on{At}(E))\longrightarrow \on{Hom}(\Ho(X,E),\Hi(X,E)).
\end{equation*}
Dually, the Petri map factors through extended Higgs bundles,
\begin{equation*}
\Ho(X,E)\ot\Ho(X,\Ev)\longrightarrow \Ho(X,\exhiggs_X(E))
\longrightarrow \Ho(X,E\ot\Ev)\, ,
\end{equation*}
where its image describes the conormal line to $\Theta_g^\circ$.

For $(X,E)\in\Theta_g^\circ$, the normalized Szeg\"o kernel gives rise
to a $\Theta_g$--twisted extended Higgs field $\sbar|_{3\Delta}$:
$\sbar_E$ is a section of $\M(E)\ot\Oo(\Theta_g)|_E$, whose value on
the diagonal $\theta(E)\on{Id}_E$ vanishes for
$(X,E)\in\Theta_g$. Since it comes as the restriction of a global kernel
$$
\sbar_E\in \Ho(\XX,E\boxt\Ev)\ot\Theta_X|_{\{E\}}\, ,
$$
it lies in the conormal bundle. Thus $\s|_{3\Delta}$ indeed has
logarithmic singularities along $\Theta_g^{\circ}$.  Again by
\propref{normalized Szego}, this gives the natural trivialization of
the twisted conormal bundle, so that the residue is $1$. We finally
conclude that the meromorphic identification of $\Omega$--torsors
$\excon_g(n)\to\Conn_g(\Theta)$ defined by $\s|_{3\Delta}\mapsto
\dlog\theta$ remains regular along $\Theta_g$, as desired. 

\subsection{Final Remarks.}
In \cite{BS}, Beilinson and Schechtman give a canonical local
description of the Atiyah bundle of the determinant line bundle, with
its Lie algebroid structure, for an arbitrary family of curves and
vector bundles. They formulate the result in the language of kernel
functions. This was the model for our definition of extended
connections. In this language, their description of the Atiyah algebra
is equivalent to a canonical, local identification of the twisted
cotangent bundle $\Conn_g(\Theta)$ with the twisted cotangent bundle
formed by the space of extended connections (equipped with a natural
symplectic form). This follows since the Atiyah algebra of a bundle
$\Ll$ is the Poisson algebra of affine--linear functions on the affine
bundle of connections on $\Ll$. The approach presented above is
global, and ignores the symplectic structure. The compatibility of our
identification with that of \cite{BS} follows as in
\secref{uniqueness} from the rigidity (\secref{rigidity}) of
$\Omega$--torsors on the moduli space of curves and bundles. In
particular the Lagrangian property of the Szeg\"o extended connections
follows, as in \corref{Lagrangian}.

\subsubsection{Relations with conformal field theory.}
The point of view of \cite{BS} is inspired by conformal field theory,
in particular the Virasoro--Kac--Moody uniformization of moduli spaces
(reviewed in \cite{book, Sugawara}).
This point of view can be expanded
to describe the role of the Szeg\"o kernel. In the abelian case, it is
well known (see \cite{Raina, KNTY}) that the Szeg\"o kernel is given
by the two--point functions of fermions twisted by the line bundle
$\Ll\in\Pg(X)\sm\Theta_X$. The conformal field theoretic explanation
for why $\dlog\theta$ comes from a kernel function on all of $\XX$ is
that it is given by a current one--point function $\langle
J(z)\rangle$ on $X$, which by operator product expansion comes as the
singular part at the diagonal of a fermion two--point function
$\langle \psi(z)\psi^*(w)\rangle$ on $\XX$ (see \cite{book} for an
algebraic development of the necessary conformal field theory). This
argument should generalize to the higher rank case, with the free
fermion replaced by $n$ free fermions. We hope to return to this in
future work.

The relation with logarithmic derivatives of theta functions only
``sees'' the restriction of the Szeg\"o kernel to $2\Delta$ or
$3\Delta$. To find similar interpretations of $\s|_{(n+1)\Delta}$ would
require an extension of the moduli of curves to ``$\mc
W$--moduli'', putatively associated to the vertex algebra ${\mc
W}({\mf sl}_{n})$ (\cite{book}), where the cotangent directions
include not only quadratic differentials but cubic, quartic, and up to
$n$--ary differentials on the curve.

\section*{Acknowledgements}
We are
very grateful to Eduard Looijenga for sending \cite{Lo} to the
second--named author in a different context. We would
like to thank Matthew Emerton, Peter Newstead and especially
Leon Takhtajan for very useful discussions. We thank the
referee for going through the paper very carefully.


\end{document}